\documentclass[12pt,oneside,openany,headsepline,halfparskip,normalheadings]{scrbook}

\usepackage[sort,numbers]{natbib}
\usepackage[a4paper,text={14cm,22cm},top=3.5cm,left=3.5cm]{geometry}
\usepackage{scrpage2,setspace}
\usepackage[latin1]{inputenc}
\usepackage[T1]{fontenc}
\usepackage{mathptmx}
\usepackage[scaled=.90]{helvet}
\usepackage{courier}
\usepackage{amssymb,amsmath,bbm}
\usepackage{graphicx}
\usepackage{epic}
\usepackage{pict2e}
\usepackage{longtable}
\usepackage{afterpage}
\usepackage{mathrsfs}
\newtheorem{theo}{Theorem}
\newtheorem{pro}{Proposition}

\newtheorem{lem}{Lemma}

\usepackage{rotating}
\usepackage{multirow}
\usepackage{setspace}
\usepackage{remreset}
\makeatletter
\@removefromreset{footnote}{chapter}
\makeatother

\usepackage{Diplomarbeit} 
\titlehead{
  \begin{center}
  {\Large }
    \end{center}
  \begin{center}    
      { \LARGE Philipp Hungerl\"{a}nder}
  \end{center}
\vspace*{1cm}
  \begin{center}
  {\LARGE }
    \end{center}
  \begin{center}
{ \LARGE Algorithms for Convex Quadratic Programming} 
  \end{center}
}
\subject{\LARGE \sc Diplomarbeit}
\author{\large Zur Erlangung des akamdemischen Grades\\
	\large Diplom-Ingenieur}
\title{}
\date{\large Technische Mathematik\\
 \large Alpen-Adria-Universit\"{a}t Klagenfurt\\
 \large Fakult\"{a}t f\"{u}r Technische Wissenschaften}
\publishers{\large\begin{flushleft} 
  Begutachter:  Univ.-Prof. Dipl.-Ing. Dr. Franz Rendl\\
  Institut f\"{u}r Mathematik \end{flushleft}
  \begin{flushright} 02/2009 \end{flushright}   }

\begin{document}
\singlespacing 
%
%

\frontmatter 
\maketitle
\chapter{Ehrenwörtliche Erklärung}

Ich erkläre ehrenwörtlich, dass ich die vorliegende wissenschaftliche Arbeit selbst\-ständig angefertigt und die mit ihr unmittelbar verbundenen Tätigkeiten selbst erbracht habe. Ich erkläre weiters, dass ich keine anderen als die angegebenen Hilfsmittel benutzt habe. Alle aus gedruckten, ungedruckten oder dem Internet im Wortlaut oder im wesentlichen Inhalt übernommenen Formulierungen und Konzepte sind gemäß den Regeln für wissenschaftliche Arbeiten zitiert und durch Fußnoten bzw. durch andere genaue Quellenangaben gekennzeichnet.\\

Die während des Arbeitsvorganges gewährte Unterstützung einschließlich signifikanter Betreuungshinweise ist vollständig angegeben.\\

Die wissenschaftliche Arbeit ist noch keiner anderen Prüfungsbehörde vorgelegt worden. Diese Arbeit wurde in gedruckter und elektronischer Form abgegeben. Ich bestätige, dass der Inhalt der digitalen Version vollständig mit dem der gedruckten Version übereinstimmt.\\

Ich bin mir bewusst, dass eine falsche Erklärung rechtliche Folgen haben wird.\\\\

(Unterschrift)  \;\;\;\;\;\;\;\;\;\;\;\;\;\;\;\;\;\;\;\;\;\;\;\;\;\;\;\;\;\;\;\;\;\;\;\;\;\;\;\;\;\;\;\;\;\;\;\;\;\;\;\;\;\;\;\;\;\;\;\;\;\;\;\;\;\;\;\;\;\;\;\;\;\;\; (Ort, Datum)


\chapter{Acknowledgements}

I am grateful to a number of people who have supported me during my studies
and in the development
of this work and it is my pleasure to mention them here.\\

    I want to thank my supervisor Franz Rendl for introducing me into the field
of optimization, for his enthusiasm about discussing mathematical
issues and for the large amount of time he devoted to my concerns. His ideas and
advice led me into active research and substantiated my diploma thesis.\\

    Furthermore I would like to thank my colleagues at the Mathematics
    Department at the Alpen-Adria-Universit\"{a}t Klagenfurt for providing me
    excellent working conditions.\\

Above all, my thanks go to my family for their interest and non-restrictive
support of my issues in all areas of life.


\tableofcontents

%
%

\mainmatter 
\chapter{Introduction}

The main interest of this diploma thesis is to describe and compare different,
practically successful solution methods for general convex quadratic problems
with arbitrary linear constraints. Therefore we first define the convex
quadratic program (QP) as 

\begin{subequations}\label{qpfd}
\begin{align}
\min_{x} \; \frac{1}{2}x^TQx+d^Tx\\ 
\text{subject to}\;\;\; Bx &= c,\\
b &\leq x \leq a,
\end{align}
\end{subequations}

where $Q$ is a positive definite $n$ $\times$ $n$ matrix, $B$ is a $m \times n$
matrix, $a$, $b$, $d \in \mathbb R^n$, and $c \in \mathbb R^m$. This problem 
has received considerable interest in the literature. We recall some of the
recent contributions. 

Solution methods like the augmented Lagrangian method, active-set methods and
interior point methods are among the most popular approaches to solve
\eqref{qpfd}, and can be traced back to the 1960's.

The so-called augmented Lagrangian method was first proposed by Hestenes
\cite{hes} and Powell \cite{pow}. More recent contributions making use of the
augmented Lagrangian idea are from Conn et al. \cite{cogo}, \cite{cogo1}, who
designed the nonlinear programming code LANCELOT, Dost\'{a}l \cite{dos},
who describes a semi-monotonic augmented Lagrangian algorithm for solving
large convex quadratic programming problems,
and Murtagh and Saunders \cite{musa}, \cite{musa1}, who developed a software
called MINOS that combines efficient
sparse-matrix techniques as in the revised simplex method with stable
quasi-Newton methods for handling the nonlinearities.

Active set methods for convex QP are the implementation of Gill and Murray
\cite{gimu1} called QPOPT, the software of Gould and Toint \cite{goto1} called
QPA, and Fletcher's code BQPD \cite{fle}. 

Another solution strategy consists in using (primal-dual) interior point
methods. Here we want to mention the recent contributions of Vanderbei
\cite{van1}, \cite{van}, \cite{vaum}, who designed the software package LOQO,
of M\'{e}sz\'{a}ros \cite{mes}, who built up the solver BPMPD, and of Gertz and
Wright \cite{gewr1}, \cite{gewr}, who developed the object-oriented software
package OOQP. Further important interior-point solvers for QP are CPLEX
\cite{ilo}, KNITRO \cite{byno}, \cite{wapl}, and MOSEK \cite{anan}.

For comparison of some of the above mentioned algorithms and methods, we refer
to the benchmarking articles of Mittelmann \cite{mit}, Dolan et
al. \cite{domo1}, \cite{domo}, and Gould and Toint \cite{goto}. 

Especially we want to mention that we developed our own contribution to solve
\eqref{qpfd} that we present in Chapter \ref{iasm}. It consists in combining
the method of 
multipliers with an infeasible active-set method. Our approach is
iterative. In each step we calculate an augmented Lagrange function. Then we
minimize this function using an infeasible active-set method that was already
successfully applied to similar problems, see the paper of Kunisch and Rendl
\cite{kure}. After this we  
update the Lagrange multiplier for the equality constraints. Finally we
try to solve \eqref{qpfd} directly, again with the infeasible active-set
method, starting from the optimal solution of the actual 
Lagrange function. Computational experience with our method indicates that
typically only few (most of the time only one) outer iterations
(multiplier-updates) and also only few (most of the time less than ten) inner
iterations (minimization of the Lagrange function and trying to solve
\eqref{qpfd} directly) are required to reach the optimal solution.

The diploma thesis is organized as follows. We close this chapter with
some notation used throughout. In Chapter \ref{prelim} we show the equivalence
of different QP problem formulations and present some important so-called
direct methods for solving equality-constrained QPs. We cover the most
important aspects for practically successful interior point methods for linear
and convex quadratic programming in Chapter \ref{cipm}. Chapter \ref{cfasm}
deals with ingredients for practically efficient feasible active set
methods. Finally Chapter \ref{iasm} provides a close description of our
Lagrangian infeasible active set method and further gives a convergence
analysis of the subalgorithms involved.


\textbf{Notation:} The following notation will be used throughout. $M :=
\{1,\ldots,m\}$ and $N := \{1,\ldots,n\}$ are two sets of integer numbers. For a
subset $A \subseteq N$ and $x \in \mathbb R^n$ we write
$x_A$ for the components of $x$ indexed by $A$, i.e. $x_A := (x_i)_{i \in
  A}$. The complement of $A$ will be denoted by $\overline{A}$. If $P$ is a
matrix and $A$ and $E$ are subsets of $N$, then $P_{A,E}$ is the submartix of
$P$, with rows indexed by $A$ and columns indexed by $E$. If $A=E$ we write
$P_A$ for $P_{A,A}$. By $P^{T}_{A,E}$ we identify the submatrix of $P^{T}$,
with rows indexed by $A$ and columns indexed by $E$. For $a,b \in \mathbb R^n$ we write $a \circ b$ to denote the vector of element-wise products, $a \circ b := (a_ib_i)_{i \in N}$.


\chapter{Preliminaries}\label{prelim}

In this chapter we show the equivalence
of different QP problem formulations in Section \ref{dpf} and
then in Section \ref{ecdm} we present some important so-called
direct methods for solving equality-constrained quadratic programs.

\section{Different Problem Formulations}\label{dpf}

The general quadratic program can be stated as

\begin{subequations}\label{gqp}
\begin{align}
\min_{x}\; \frac {1}{2}x^{\top}Qx &+ x^{\top}d \label{gqpof}\\
\text{subject to  } a_{i}^{\top}x &= c_{i}, \;\;\; i \in
\epsilon, \label{gqpec}\\ 
a_{i}^{\top}x &\leq c_{i}, \;\;\; i \in \iota, \label{gqpic}
\end{align}
\end{subequations}

where $Q$ is a symmetric $n \times n$ matrix, $\epsilon$ and $\iota$ are
finite sets of indices, and $d$, $x$ and $\{a_{i}\},\;\; i \in \epsilon \cup
\iota ,$ are vectors in $\mathbb R^{n}$. If the Hessian matrix $Q$ is positive
definite, we say that \eqref{gqp} is a strictly convex QP,
and in this case the problem is often similar in difficulty to a linear
program. Nonconvex QPs, in which $Q$ is an indefinite matrix, can be more
challenging because they can have several stationary points and local minima. 

We can convert the inequality constraints in the above formulation of a QP by
introducing a vector of slack variables $z$ and writing

\begin{subequations}\label{gqpsv}
\begin{align}
\min_{x} \; \frac {1}{2}x^{\top}Qx &+ x^{\top}d \\
\text{subject to  } a_{i}^{\top}x &= c_{i}, \;\;\; i \in \epsilon, \\
a_{i}^{\top}x + z &= c_{i}, \;\;\; i \in \iota,\\
z_{i} &\geq 0, \;\;\; i \in \iota.
\end{align}
\end{subequations}

We can further transform this formulation by splitting $x$ into its
nonnegative and nonpositive parts, $\;\;x = x^{+} - x^{-}\;\;$, where $x^{+} =
\max(x,0) \geq 0$ and $x^{-}= \max(-x,0) \geq 0$. The problem \eqref{gqpsv}
can now be written as

\begin{align*}
\min_{(x^{+},x^{-},z)} \; \frac {1}{2}\begin{pmatrix}x^{+} \\ x^{-} \\
  z \end{pmatrix}^{\top}&\begin{pmatrix}Q & 0 & 0 \\ 0 & Q & 0 \\ 0 & 0 &
  0 \end{pmatrix}\begin{pmatrix}x^{+} \\ x^{-} \\ z \end{pmatrix}
+ \begin{pmatrix}x^{+} \\ x^{-} \\ z \end{pmatrix}^{\top}\begin{pmatrix} d \\
  -d \\ 0 \end{pmatrix} \\ 
\text{subject to  } \begin{pmatrix}a_{i} \\ -a_{i} \\
  0 \end{pmatrix}^{\top}\begin{pmatrix}x^{+} \\ x^{-} \\
  z \end{pmatrix} &= c_{i}, \;\;\; i \in \epsilon \\ 
\begin{pmatrix}a_{i} \\ -a_{i} \\ 1 \end{pmatrix}^{\top} \begin{pmatrix}x^{+}
  \\ x^{-} \\ z \end{pmatrix} &= c_{i}, \;\;\; i \in \iota\\
\begin{pmatrix}x^{+} \\ x^{-} \\ z \end{pmatrix} &\geq 0.
\end{align*}

Now setting 

\begin{align*}
\begin{pmatrix}x^{+} \\ x^{-} \\ z \end{pmatrix} = \overline{x},\\
\begin{pmatrix}Q & 0 & 0 \\ 0 & Q & 0 \\ 0 & 0 & 0 \end{pmatrix} =
\overline{Q}, \\ 
\begin{pmatrix} d \\ -d \\ 0 \end{pmatrix} = \overline{d}, \\
\begin{pmatrix} A_{\epsilon} & -A_{\epsilon} & 0_{n \times k} \\ A_{\iota} &
  -A_{\iota} & I_{n \times m-k} \end{pmatrix} = B,
\end{align*}

where

\begin{align*}
A_{\epsilon} &= [a_{i}]_{i \in \epsilon},\\
A_{\iota} &= [a_{i}]_{i \in \iota},\\
c &= [c_{i}]_{i \in \epsilon \cup \iota}, \\
k &= |\epsilon|,\\
m &= |\iota| + |\epsilon|,
\end{align*}

we obtain

\begin{subequations}\label{qpsf}
\begin{align}
\min_{\overline{x}} \; \frac {1}{2}\overline{x}^{\top}\overline{Q}\overline{x}
&+ \overline{x}^{\top}\overline{d} \\ 
\text{subject to} \;\;B \overline{x} &= c,\\
\overline{x} &\geq 0.
\end{align}
\end{subequations}

Hence, we showed that \eqref{qpsf} is equivalent to \eqref{gqp} and it depends
on the considered algorithm what representation of the quadratic problem is
preferable. Furthermore we want to mention that we can also convert inequality
constraints of the form $x \leq a$ or $Ax \geq c$ to equality constraints by
adding or subtracting slack variables: 

\begin{align*}
x \leq a \;\; \Leftrightarrow \;\; x + w = a, \;\; w \geq 0, \\
Ax \geq c \;\; \Leftrightarrow \;\; Ax - w = c, \;\; w \geq 0.
\end{align*}


\section{Solution Methods For Equality-Constrained QPs} \label{ecdm}

In this section we consider direct solution methods for quadratic programs in
which only equality constraints are present. We define them as follows:

\begin{subequations} \label{ecqp}
\begin{align}
\min_{x} \frac {1}{2}x^{\top}Qx + x^{\top}d \\
\text{subject to} \;\;\; Ax = c,
\end{align}
\end{subequations}

where $Q$ is a positive-definite $n \times n$ matrix, $A$ is a $m \times n$
non-singular matrix, $d$ is a vector in $\mathbb R^{n}$ and $c$ is a vector in
$\mathbb R^{m}$. 

The KKT conditions for this problem are 

\begin{align}
 K \begin{bmatrix} -p\\ \lambda^{*} \end{bmatrix} =
\begin{bmatrix} g\\ h \end{bmatrix} \label{ecqpkkt}
\end{align}

where 

\begin{align}
K = \begin{bmatrix} Q & A^{\top}\\ A &  0 \label{ecqpkktm}\end{bmatrix}
\end{align}

and

\begin{align*}
g &= c + Qx, \\
h &= Ax - b, \\
p &= x^{*} - x.
\end{align*}

These problems appear often as subproblems in algorithms that solve general
QPs with inequality constraints (see, for example, the subproblems for feasible
active-set methods 
described in Section \ref{asmsp}) and therefore it is very important to find
ways to solve them efficiently. In addition to the direct solution methods
described in this section, there also exist iterative solution methods like
the conjugate gradient method applied to the reduced system and the projected
conjugate gradient method. For a further discussion of these iterative
methods see, for example, Nocedal and Wright \cite[Section 16.3]{nowr}, Conn,
Gould, and Toint \cite{cogo1} and Burke and Mor\'{e} \cite{bumo}.

\subsection{Factoring the full KKT system}

One option for solving \eqref{ecqpkkt} is the use of a triangular
factorization of $K$ and then make backward and forward substitution. 
To discuss this option we need some theoretical knowledge about the
definiteness of $K$. Therefore let us give a result that states that the KKT
matrix $K$ is always indefinite. We define

\begin{align*}
\text{inertia}(S) \stackrel{\mathrm{def}}= (n_{+},n_{-},n_{0})
\end{align*}

where $n_{+}$ is the number of positive eigenvalues of $S$, $n_{-}$ denotes the
number of negative eigenvalues of $S$ and $n_{0}$ gives the number of zero
eigenvalues of $S$. Now we can state a result that characterizes the inertia of $K$.

\begin{theo} Let K be given by \eqref{ecqpkktm}, and suppose that A has rank
  m, Then

\begin{align*}
\text{inertia}(K) = \text{inertia}(Z^{\top}QZ) + (m,m,0),
\end{align*}

where $Z$ is an $n \times (n-m)$ matrix whose columns are a basis of the null
space of A. That is, $Z$ has full rank and satisfies $AZ = 0$. Therefore, if $Z^{\top}QZ$ is positive definite, \text{inertia}(K) = (n,m,0).
\label{ink}
\end{theo}

The proof of this result is given in Forsgren and Gill \cite[Lemma 4.1]{fogi}
or Gould \cite[Lemma 3.4]{gou}, for example.

Because of indefiniteness of $K$, we cannot use the Cholesky
factorization to solve \eqref{ecqpkkt}. The use of Gaussian elimination has the
disadvantage that it 
ignores symmetry. Therefore the most effective approach is to use a symmetric
indefinite factorization\footnote{The computational cost of a symmetric
  factorization is typically about half the cost of Gaussian elimination}
which has the form 

\begin{align*}
\bar{P}^{\top}S\bar{P} = \bar{L}\bar{B}\bar{L}^{\top}  
\end{align*}

where $S$ is a general symmetric matrix, $\bar{P}$ is a permutation
matrix, $\bar{L}$ is a  unit lower triangular matrix and $\bar{B}$
is a block-diagonal matrix with either 1 
$\times$ 1 or 2 $\times$ 2 blocks. We use the symmetric permutations defined
by $\bar{P}$ to improve the numerical stability of the computation and,
if $S$ is sparse, to maintain sparsity. 

Now to solve \eqref{ecqpkkt}, we first compute a factorization of KKT matrix
$K$:

\begin{align}
P^{\top}KP = LBL^{\top}, \label{fac}
\end{align}

and then use the calculated factors in the following way to arrive the
solution:

\begin{align*}
\text{solve} \;\; Lz &= P^{\top}\begin{bmatrix} g \\ h \end{bmatrix} \;\;
\text{to obtain} \;\; z; \\
\text{solve} \;\; B\hat{z} &= z \;\;\;\;\;\;\;\;\;\;\;\;\;\; \text{to obtain}\;\; \hat{z}; \\
\text{solve} \;\; L^{\top}\bar{z} &= \hat{z} \;\;\;\;\;\;\;\;\;\;\;\;\;\; \text{to obtain}\;\;
\bar{z}; \\ 
\text{set} \;\; \begin{bmatrix} -p \\ \lambda^{*} \end{bmatrix} &= P \bar{z}.
\end{align*}

The by far most expensive operation in this approach is the performance of the
factorization \eqref{fac}. This factoring of the KKT matrix $K$ is quite
effective for many problems. It may be expensive,
however, if $K$ is sparse and the heuristics for choosing $P$ are not able to
maintain this sparsity in $L$ and therefore $L$ becomes dense.

\subsection{The Schur-complement method}

We assumed for the equality-constrained QP \eqref{ecqp} that $Q$ is positive
definite. Therefore we can multiply the first equation in \eqref{ecqpkkt} by
$AQ^{-1}$ and then subtract the second equation to get the following equation
in $\lambda^{*}$ alone:

\begin{align*}
(AQ^{-1}A^{\top})\lambda^{*} = (AQ^{-1}g - h).
\end{align*}

As $AQ^{-1}A^{\top}$, the so called Schur complement of $Q$, is also positive
definite (because we assumed that $A$ has full rank), we can calculate
$\lambda^{*}$ as:

\begin{align*}
\lambda^{*} = (AQ^{-1}A^{\top})^{-1}(AQ^{-1}g - h),
\end{align*}

and then obtain $p$ from the first equation of \eqref{ecqpkkt}:

\begin{align*}
p = Q^{-1}(A^{\top}\lambda^{*} - g).
\end{align*}

Using the Schur-Complement method we need to invert $Q$, as well as to compute
a factorization of the $m \times m$ matrix $AQ^{-1}A^{\top}$. Therefore, the
method is most effective if $Q$ is well conditioned and easy to invert or if
$Q^{-1}$ is known explicitly through a quasi-Newton updating formular or if
the number of constraints $m$ is small.\\

\subsection{The null-space method}

The null-space method does not require nonsingularity of $Q$ but only full
rank of $A$ and positive definiteness of $Z^{\top}QZ$, where $Z$ is the
null-space basis matrix.

Let us partition the vector $p$ in \eqref{ecqpkkt} into two components, so
that:

\begin{align}
p = Yp_{y} + Zp_{z}, \label{psp} 
\end{align}

where $Z$ is $n \times (n-m)$, Y is $n \times m$, $p_{y}$ is a vector in
$\mathbb R^{m}$ and $p_{z}$ is a vector in $\mathbb R^{n-m}$. 

Thereby we choose $Y$ and $Z$ with the following properties:

\begin{align*}
\begin{bmatrix} Y | Z \end{bmatrix} \in \mathbb R^{n \times n} \;\; \text{is
  nonsingular,} \;\; AZ = 0.
\end{align*}

Since $A$ has full rank, so does $A[Y|Z] = [AY|0]$ and therefore $AY$ is
nonsingular and has rank $m$.

Now we substitute $p$ with the help of \eqref{psp} in $Ap = -h$, which gives

\begin{align*}
(AY)p_{y} = -h.
\end{align*}

We can make $p_{y}$ explicit, as $AY$ is nonsingular:

\begin{align}
p_{y} = -(AY)^{-1}h. \label{defpy}
\end{align}

To determine $p_{z}$ we use the first equation of \eqref{ecqpkkt} to obtain

\begin{align*}
-QYp_{y} -QZp_{z} + A^{\top}\lambda^{*} = g,
\end{align*}

and then multiply it by $Z^{\top}$:

\begin{align}
(Z^{\top}QZ)p_{z} = - Z^{\top}QYp_{y} - Z^{\top}g. \label{Zrs}
\end{align}

To calculate $p_{z}$ from this equation, we can use, for example, a Cholesky
factorization of $Z^{\top}QZ$. After that we can compute the total step $p$ by
using \eqref{psp}. Finally we can obtain $\lambda^{*}$ by multiplying the
first equation of \eqref{ecqpkkt} by $Y^{\top}$

\begin{align*}
(AY)^{\top}\lambda^{*} = Y^{\top}(g + Qp), \label{lamopt}
\end{align*}

and then solving this equation for $\lambda^{*}$.

The main computational effort of the null-space method lies in the
determination of the not uniquely defined matrix $Z$. If we choose $Z$ to have
orthonormal columns\footnote{For a orthonormal $Z$ the corresponding $Y$ can be
  calculated by a QR factorization of $A^{\top}$, for details see Section
  \ref{udf}}, then the conditioning of $Z^{\top}QZ$ is at least as good as
that of $Q$ itself, but a orthonormal $Z$ is often expensive to compute
(especially if $A$ is sparse). On the other hand if we choose $Z$ in a
different, computationally cheaper way, the reduced system \eqref{Zrs} may
become ill 
conditioned. Therefore the null-space method is preferable compared with
the Schur-complement method when it is more expensive to invert $Q$ and
compute factors 
of $AQ^{-1}A^{\top}$ than to compute $Z$ and factors of $Z^{\top}QZ$ and
$AY$. This is most of the time the case if the number of equality constraints
$m$ is large and therefore the matrices $Z$ and $Z^{\top}QZ$ have low
dimensions.


\chapter{Interior Point Methods}\label{cipm}

This chapter is devoted to the description of practically successful interior
point methods for linear and convex quadratic programming. In Section
\ref{ashr} we mention some basis information about the exciting historical
development of the interior point methods as first real competitor of the
simplex method. Section \ref{slp} is used to present the central components of
interior point methods on the basis of the simple linear programming
framework. After that we show in Section \ref{etcqp} that the generalisation 
of interior point methods to QPs is a natural and easy one, especially if we
compare it with the serious differences between the simplex method and active
set methods for QPs.

\section{A short historical review}\label{ashr}

Starting with the seminal paper of Karmarkar \cite{kar} in 1984, interior
point methods in mathematical programming have been the most 
important research area in optimization since the development of the
simplex method for linear programming. Interior point methods have strongly
influenced mathematical programming theory, practice and computation. For
example linear programming is no longer synonymous with the simplex method,
and linear programming is shown as a special case of nonlinear programming due
to these developments. 

On the theoretical side, permanent research led to better computational
complexity bounds for linear programming, quadratic programming, linear
complementarity problems, semi-definite programming and some classes of convex
programming problems. On the computational side, the performance of tools for
linear and nonlinear programming improved greatly, as the sudden appearance of
credible competition for the active set methods initiated significant
improvements in implementations.

Interior-point methods arose from the search for algorithms with
better theoretical properties than the simplex method. As Klee and Minty
\cite{klmi} showed, the simplex method can be inefficient on certain
pathological problems. Roughly speaking, the time required to solve a linear
program may be exponential in the size of the problem, as measured by the
number of unknowns and the amount of storage needed for the problem data. For
almost all practical problems, the simplex method is much more efficient than
this bound would suggest, but its poor worst-case complexity motivated the
development of new algorithms with better guaranteed performance. The first
such method was the ellipsoid method, proposed by Khachiyan \cite{kha}, which
finds solution in time that is at worst polynomial in the problem
size. Unfortunately, this method approaches its worst-case bound on all
problems and is not competitive with the simplex method in practice.

Karmarkar's projective algorithm \cite{kar}, announced in 1984, also has the
polynomial complexity property, but it came with the added attraction of good
practical behavior. The initial claims of excellent performance on large
linear programs were never fully borne out, but the announcement prompted a
great deal of research activity which gave rise to many new methods.

In the first years after Karmarkar's initial paper, research in linear
programming was concentrated on finding algorithms that worked with the primal
problem, but had better complexity bounds or were easier to implement than the
original method. A next crucial step was done by Megiddo \cite{meg} in 1987,
when he described a framework for primal-dual algorithms. To take into account
the primal and the dual problem proved to be extraordinarily productive. The
primal-dual viewpoint led to new algorithms with best practical and also
interesting theoretical properties. Furthermore it formed the basis for
transparent extensions to convex programming and linear complementarity. The
basis algorithm for most current practical linear programming software was
described by Mehrotra in 1989 \cite{meh}.  

Some years later, Nesterov and
Nemirovskii published their theory of self-concordant functions \cite{nene}
which was the main tool to extend algorithms for linear programming based on
the primal log-barrier function to more general classes of convex problems
like semi-definite programming and second-order cone programming. Later on,
Nesterov and Todd \cite{neto,neto1} did further extending work along these
lines. Interior point methods have also been frequently used in such areas as
control theory, structural optimization, combinatorial and integer programming
and linear algebra for different decomposition methods.

In the next sections, we will concentrate on central trajectory methods using
the primal-dual framework, because these algorithms have the best practical
features in the class of interior point methods. Furthermore we will
concentrate on linear and convex quadratic programming. Readers interested
also in affine scaling or potential reduction methods or in algorithms using
only the primal or only the dual variables or in further topics like linear
complementarity problems, semi-definite programming, self-duality or and
theoretical run-time properties are referred to three survey articles of
Forsgren et al. \cite{fogi1}, Freund and Mizuno \cite{frmi} and Potra and
Wright \cite{powr} and two comprehensive books of Wright \cite{wri} and Ye
\cite{ye} about interior point methods.


\section{Linear Programming}\label{slp}

In this section we present the central components of interior point methods on
the basis of the simple linear programming framework. We consider the linear
programming problem in standard form 

\begin{subequations}\label{linpropri}
\begin{align}
\min c^{\top}x \\
\text{subject to } Ax &= b, \\
x &\geq 0, 
\end{align}
\end{subequations}

where $c$ and $x$ are vectors in $\mathbb R^{n}$, $b$ is a vector in $\mathbb
R^{m}$, and $A$ is an $m \times n$ matrix with full row rank. The dual problem
for \eqref{linpropri} is

\begin{subequations}\label{linprodua}
\begin{align}
\max b^{\top}\lambda\\
\text{subject to} \;\; A^{\top}\lambda + s &= c, \\
s &\geq 0, 
\end{align}
\end{subequations}

where $\lambda$ is a vector in $\mathbb R^{m}$ and s is a vector in $\mathbb
R^{n}$. 

\subsection{The KKT system and Newton's method}

Solutions of \eqref{linpropri} together with \eqref{linprodua} are
characterized by the KKT conditions:

\begin{subequations}\label{KKTIPl}
\begin{align}
A^{\top}\lambda + s &= c, \label{IPlinopt}\\
Ax &= b, \label{IPlincon}\\
x_{i}s_{i} &= 0,\;\;\;\; i = 1,2,\ldots,n, \label{IPlincomp}\\
x &\geq 0, \label{IPlinxb}\\
s &\geq 0. \label{IPlinsb}
\end{align}
\end{subequations}

Primal-dual methods find solutions $(x^{*},\lambda^{*},s^{*})$ of this system
by applying variants of Newton's method to the three equalities
\eqref{IPlinopt} - \eqref{IPlincomp} and modifying the search directions and
step lengths so that the inequalities \eqref{IPlinxb} and \eqref{IPlinsb} are
satisfied strictly in 
every iteration. The equations \eqref{IPlinopt} and \eqref{IPlincon} are
linear and \eqref{IPlincomp} is only mildly nonlinear. So these three
equations are not difficult to solve by themselves. However, the problem
becomes much more difficult when we add the nonnegativity requirements
\eqref{IPlinxb} and \eqref{IPlinsb}, which give rise to all the complications
in the design and analysis of interior-point methods.

To derive primal-dual interior-point methods we restate the first three
equations \eqref{IPlinopt} - \eqref{IPlincomp} of the above KKT-system
in a slightly different form by means of a mapping $F$ from $\mathbb R^{2n+m}$
to $\mathbb R^{2n+m}$:

\begin{align}
F(x,\lambda,s) = \begin{bmatrix} A^{\top}\lambda + s - c \\ Ax - b \\
  XSe \end{bmatrix} &= 0 \label{IPlinF}
\end{align}

where

\begin{align*}
X &= diag(x_{1},\ldots,x_{n}),\\
S &= diag(s_{1},\ldots,s_{n}),
\end{align*}

and $e = (1,\ldots,1)^{\top}$. Primal-dual methods generate iterates
$(x^{k},\lambda^{k},s^{k})$ that satisfy the bounds \eqref{IPlinxb} and
\eqref{IPlinsb} strictly. This property is the origin of the term
interior-point. By respecting these bounds, the method avoids solutions, that
satisfy $F(x,\lambda,s)$ = 0 but not \eqref{IPlinxb} or \eqref{IPlinsb}. These
so-called spurious solutions abound and do not provide useful information
about solutions of \eqref{linpropri} or \eqref{linprodua}, so it makes sense
to exclude them altogether from the region of search.

\subsection{The duality measure and the centering parameter}

Like most iterative algorithms in optimization, primal-dual interior-point
methods have two basic ingredients; a procedure for determining the step and a
measure of the desirability of each point in the search space. An important
component of the measure of desirability is the average value of the pairwise
products $x_{i}s_{i},\;\; i = 1,\ldots,n$, which are all positive when $x>0$ and
$s>0$. This quantity is known as the duality measure and is defined as
follows:

\begin{align}
\mu = \frac{1}{n}\sum_{i=1}^{n}x_{i}s_{i} = \frac{x^{\top}s}{n}. \label{linpromu}
\end{align}

The procedure for determining the search direction has its origins in Newton's
method for the nonlinear equations \eqref{IPlinF}. Newton's method forms a
linear model of \eqref{IPlinF} around the current point and obtains the search
direction $(\Delta x,\Delta \lambda, \Delta s)$ by solving the following
system of linear equations:

\begin{align}
\begin{bmatrix} 0 & A^{\top} & I \\ A & 0 & 0 \\ S & 0 & X \end{bmatrix}
\begin{bmatrix} \Delta x \\ \Delta \lambda \\ \Delta s \end{bmatrix} = 
\begin{bmatrix} -r_{c} \\ -r_{b} \\ -XSe \end{bmatrix} \label{linpune}
\end{align}

where

\begin{align*}
r_{b} &= Ax - b,\\ 
r_{c} &= A^{\top}\lambda + s - c.
\end{align*}

Usually, a full step along this direction would violate the bounds, so we
perform a line search along the Newton direction and define the new iterate as

\begin{align*}
(x^{+}, \lambda^{+} ,s^{+}) = (x,\lambda,s) + \alpha (\Delta x,\Delta \lambda,
\Delta s), 
\end{align*}

for some line search parameter $\alpha \in (0,1]$. We often can take only a
small step along this direction ($\alpha << 1$) before violating the
conditions $x>0$ and $s>0$. Hence, the pure Newton direction \eqref{linpune},
sometimes known as the affine scaling direction, often does not allow us to
make much progress towards a solution.

Most primal-dual methods use a less aggressive Newton direction, one that does
not aim directly for a solution of \eqref{IPlinF} but rather for a point whose
pairwise products $x_{i}s_{i}, \;\; i \in 1,\ldots,n,$ are reduced to a lower average value - not all
the way to zero. So we take a Newton step toward a point for which $x_{i}s_{i}
= \sigma \mu$, where $\mu$ is defined by \eqref{linpromu} and $\sigma \in
[0,1]$ is the reduction factor that we wish to achieve in the duality measure
on this step. The modified step equation is then

\begin{align}
\begin{bmatrix} 0 & A^{\top} & I \\ A & 0 & 0 \\ S & 0 & X \end{bmatrix}
\begin{bmatrix} \Delta x \\ \Delta \lambda \\ \Delta s \end{bmatrix} = 
\begin{bmatrix} -r_{c} \\ -r_{b} \\ -XSe + \sigma \mu e \end{bmatrix}. \label{linsine}
\end{align}

When $\sigma > 0$, it usually is possible to take a longer step $\alpha$ along
the direction defined by \eqref{linsine} before violating the
bounds. Therefore $\sigma$ is called the centering parameter.

The concrete choices of $\sigma$ and $\alpha$ are crucial to the performance
of interior-point methods. Therefore techniques for controlling these
parameters, directly and indirectly, give rise to a wide variety of methods
with diverse properties.

Although software for implementing interior point methods does usually not
start from a feasible point $(x^{0},\lambda^{0},s^{0})$ that 
fulfills:\footnote{For a feasible starting point we have $r_{b}$ = $r_{c}$ = 0
for all iterations of the so-called feasible interior point method.} 

\begin{align*}
Ax^{0} = b,\\
A^{\top}\lambda^{0} + s^{0} = c,
\end{align*} 

most of the historical development of theory and algorithms assumed that these
conditions are satisfied. Assuming this, a comprehensive convergence analysis
can be 
presented in just a few pages using only basic mathematical tools and concepts
(see, for example, Wright \cite[Chapter 5]{wri} or Nocedal and Wright
\cite[Section 14.1]{nowr}). Analysis of the infeasible case follows the same
principles, but is considerably more complicated in the details. 

Practical implementations of interior-point algorithms work with an infeasible
starting point and infeasible iterations, maintain strict positivity of $x$
and $s$ throughout and take at each iteration a Newton-like step involving a
centering component. Several aspects of 'theoretical` algorithms are typically
ignored, while several enhancements are added that have a significant effect
on practical performance. Next we describe the algorithmic enhancements that
are found in a typical implementation of an infeasible-interior-point method
(for further details consult the paper of Mehrotra \cite{meh}).

\subsection{Corrector and centering steps}

A key feature of practical algorithms is their use of corrector steps that
compensate for the linearization error made by the Newton affine-scaling
step in modeling equation \eqref{IPlincomp}. Consider the affine-scaling
direction $(\Delta x^{aff}, \Delta \lambda^{aff}, \Delta s^{aff})$ defined by

\begin{align}
\begin{bmatrix} 0 & A^{\top} & I \\ A & 0 & 0 \\ S & 0 & X \end{bmatrix}
\begin{bmatrix} \Delta x^{aff} \\ \Delta \lambda^{aff} \\ \Delta s^{aff} \end{bmatrix} = 
\begin{bmatrix} -r_{c} \\ -r_{b} \\ -XSe \end{bmatrix}. \label{linpupre}
\end{align}

If we take a full step in this direction, we obtain

\begin{align*}
(x_{i} + \Delta x_{i}^{aff})(s_{i} + \Delta s_{i}^{aff}) = x_{i}s_{i} + x_{i}
\Delta s_{i}^{aff} + s_{i} \Delta x_{i}^{aff} + \Delta x_{i}^{aff} \Delta
s_{i}^{aff} = \Delta x_{i}^{aff} \Delta s_{i}^{aff}.
\end{align*}

That is, the updated value of $x_{i}s_{i}$ is $\Delta x_{i}^{aff} \Delta
s_{i}^{aff}$ rather than the ideal value 0. We can solve the following system
to obtain a step $(\Delta x^{cor}, \Delta \lambda^{cor}, \Delta s^{cor})$ that
attempts to correct for this deviation form the ideal:

\begin{align}
\begin{bmatrix} 0 & A^{\top} & I \\ A & 0 & 0 \\ S & 0 & X \end{bmatrix}
\begin{bmatrix} \Delta x^{cor} \\ \Delta \lambda^{cor} \\ \Delta s^{cor} \end{bmatrix} = 
\begin{bmatrix} 0 \\ 0 \\ -\Delta X^{aff} \Delta S^{aff} e \end{bmatrix} \label{linuncor}
\end{align}

In many cases, the combined step $(\Delta x^{aff}, \Delta \lambda^{aff},
\Delta s^{aff}) + (\Delta x^{cor}, \Delta \lambda^{cor}, \Delta s^{cor})$ does
a better job of reducing the duality measure than does the affine-scaling step
alone.

A further important ingredient for a good practical algorithm is the use of
centering steps, with an adaptive choice of the centering parameter
$\sigma_{k}$. Thereby the affine-scaling step can be used as the basis of a
successful heuristic for choosing $\sigma_{k}$.

If the affine-scaling step reduces the duality measure significantly, there is
not much need for centering, so a smaller value of $\sigma_{k}$ is
appropriate. Conversely, if not much progress can be made along this direction
before reaching the boundary of the nonnegative orthant, a larger value of
$\sigma_{k}$ will ensure that the next iterate is more centered, so a longer
step will be possible form this next point. The following scheme calculates
the maximum allowable step lengths along the affine-scaling direction:

\begin{subequations}\label{alphaffpridua}
\begin{align}
\alpha_{aff}^{pri} \stackrel{\mathrm{def}}= \min(1, \min_{i:\Delta
  x_{i}^{aff}<0}- \frac {x_{i}}{\Delta x_{i}^{aff}}), \\
\alpha_{aff}^{dual} \stackrel{\mathrm{def}}= \min(1, \min_{i:\Delta
  s_{i}^{aff} < 0} - \frac {s_{i}}{\Delta s_{i}^{aff}}). 
\end{align}
\end{subequations}

Now we can define $\mu_{aff}$ to be the value of $\mu$ that would be obtained
by using these step lengths:

\begin{align}
\mu_{aff} = (x + \alpha_{aff}^{pri} \Delta x^{aff})^{\top}(s +
\alpha_{aff}^{dual} \Delta s^{aff})/n. \label{muaff}
\end{align}

The centering parameter $\sigma$ is chosen according to the following
heuristic:\footnote{\eqref{sigheu} has no solid analytical justification, but
  appears to work well in practice}

\begin{align}
\sigma = (\frac{\mu_{aff}}{\mu})^{3}. \label{sigheu}
\end{align}

To summarize, computation of the search direction requires the solution of two
linear systems. First \eqref{linpupre} is solved to obtain the affine-scaling
direction, also known as the predictor step. This step is used to define the
right-hand side for the corrector step and to calculate the centering
parameter from \eqref{alphaffpridua} - \eqref{sigheu}. Second, the search
direction is calculated solving

\begin{align}
\begin{bmatrix} 0 & A^{\top} & I \\ A & 0 & 0 \\ S & 0 & X \end{bmatrix}
\begin{bmatrix} \Delta x \\ \Delta \lambda \\ \Delta s \end{bmatrix} = 
\begin{bmatrix} -r_{c} \\ -r_{b} \\ -XSe - \Delta X^{aff} \Delta S^{aff} e +
  \sigma \mu e \end{bmatrix}. \label{linfindi}
\end{align}

Note that the predictor, corrector and centering contributions have been
aggregated on the right-hand side of this system. The coefficient matrix in
both linear systems \eqref{linpupre} and \eqref{linfindi} is the same. Thus,
the factorization of the matrix needs to be computed only once and the
marginal cost of solving the second system is relatively small.

\subsection{Finding an optimal step length}

Practical implementations typically calculate the maximum step lengths that
can be taken in the $x$ and $s$ variables without violating
nonnegativity separately:

\begin{subequations}
\begin{align*}
\alpha_{k,max}^{pri} &\stackrel{\mathrm{def}} = \min(1, \min_{i:\Delta
  x_{i}^{k}<0}- \frac {x_{i}}{\Delta x_{i}^{k}}), \\
\alpha^{dual}_{k,max} &\stackrel{\mathrm{def}} = \min(1, \min_{i:\Delta
  s_{i}^{k} < 0} - \frac {s_{i}}{\Delta s_{i}^{k}}). 
\end{align*}
\end{subequations}

and then take a step length of slightly less than this maximum:

\begin{subequations}\label{alphpridua}
\begin{align}
\alpha_{k}^{pri} &= \eta_{k}\alpha_{k,max}^{pri},\\
\alpha_{k}^{dual} &= \eta_{k}\alpha_{k,max}^{dual},
\end{align}
\end{subequations}

where $\eta_{k} \in [0.9,1.0]$ is chosen to accelerate the asymptotic
convergence. Therefore we want that $\eta_{k} \rightarrow 1$ as
the iterates approach the primal-dual solution. Then we obtain a new iterate
by setting 

\begin{align*}
x^{k+1} &= x^{k} + \alpha_{k}^{pri} \Delta x^{k}, \\
(\lambda^{k+1},s^{k+1}) &= (\lambda^{k},s^{k}) + \alpha_{k}^{dual}(\Delta
\lambda^{k}, \Delta s^{k}).
\end{align*}

As the step $(\Delta x^{k}, \Delta \lambda^{k}, \Delta s^{k})$ adjusts the
infeasibility in the KKT conditions

\begin{align*}
A \Delta x^{k} &= -r_{b}^{k}, \\
A^{\top} \Delta \lambda^{k} + \Delta s^{k} &= -r_{c}^{k}
\end{align*}

we have for the infeasibilities at the new iterate $k+1$

\begin{align*}
r_{b}^{k+1} &= (1-\alpha_{k}^{pri})r_{b}^{k}, \\
r_{c}^{k+1} &= (1-\alpha_{k}^{dual})r_{c}^{k}.
\end{align*}

\subsection{Choosing a starting point}

Choice of the starting point is an important practical issue with a significant
effect on the robustness of the algorithm. A poor choice
$(x^{0},\lambda^{0},s^{0})$ satisfying only \eqref{IPlinxb} and
\eqref{IPlinsb} often leads to failure in convergence. We describe here a
heuristic that finds a starting point that satisfies \eqref{IPlinopt} and
\eqref{IPlincon} reasonably well, while maintaining \eqref{IPlinxb} and
\eqref{IPlinsb} and additionally avoiding too large values of these
components.

First we find a vector $\tilde{x}$ of minimum norm satisfying \eqref{IPlincon}:

\begin{align*}
\min_{x} \frac {1}{2} x^{\top}x \\
\text{ subject to } Ax = b,
\end{align*}

and a vector $(\tilde{\lambda},\tilde{s})$ satisfying 
\eqref{IPlinopt} such that $\tilde{s}$ has minimum norm:

\begin{align*}
\min_{(\lambda,s)} \frac {1}{2} s^{\top}s \\
\text{ subject to } A^{\top}\lambda + s = c.
\end{align*}

The optimal values $(\tilde{x}, \tilde{\lambda},\tilde{s})$
can be written explicitly as follows:

\begin{subequations}\label{IPspopt} 
\begin{align}
\tilde{x} &= A^{\top}(AA^{\top})^{-1}b,\\
\tilde{\lambda} &= (AA^{\top})^{-1}Ac, \\
\tilde{s} &= c - A^{\top}\tilde{\lambda}.
\end{align}
\end{subequations}

In general, $\tilde{x}$ and $\tilde{s}$ will have nonpositive components, so
are not suitable for use as a starting point. Therefore define

\begin{align*}
\delta_{x} &= \max(- \frac {3}{2} \min_{i} \tilde{x}_{i}, 0), \\
\delta_{s} &= \max(- \frac {3}{2} \min_{i} \tilde{s}_{i}, 0),
\end{align*}

and adjust the $\tilde{x}$ and $\tilde{s}$ vectors so that they fulfill
\eqref{IPlinxb} and \eqref{IPlinsb}:

\begin{align*}
\hat{x} &= \tilde{x} + \delta_{x}e,\\
\hat{s} &= \tilde{s} + \delta_{s}e.
\end{align*}

To ensure that our starting points are not too close to zero and not too
dissimilar, we define them finally as:

\begin{subequations}\label{IPfsp}
\begin{align}
x^{0} &= \hat{x} + \hat{\delta_{x}}e, \\
\lambda^{0} &= \tilde{\lambda}, \\
s^{0} &= \hat{s} + \hat{\delta_{s}}e. 
\end{align}
\end{subequations}

where

\begin{align*}
\hat{\delta_{x}} &= \frac {1}{2} \frac{\hat{x}^{\top}\hat{s}}{e^{\top}\hat{s}},\\
\hat{\delta_{s}} &= \frac{1}{2} \frac {\hat{x}^{\top}\hat{s}}{e^{\top}\hat{x}}
\end{align*}

The computational cost of finding $(x^{0},\lambda^{0},s^{0})$ by this scheme
is about the same as one step of the primal-dual method.\label{lpsp}

\subsection{A practical primal-dual method}

Finally we put together the different, above mentioned ingredients for a
practically successful algorithm in Table \ref{lpppca} below.\\

\begin{table}[h]
\small 
\begin{tabular}{l} \hline
\multicolumn{1}{c}{Practical Predictor-Corrector Algorithm} \\\hline
Calculate $(x^{0},\lambda^{0},s^{0})$ using \eqref{IPspopt} - \eqref{IPfsp};
\\\hline 
k = 0\\
\textbf{repeat}\\
Set $(x,\lambda,s)$ = $(x^{k},\lambda^{k},s^{k})$ and calculate $(\Delta
x^{aff}, \Delta \lambda^{aff}, \Delta s^{aff})$ by solving \eqref{linpupre};\\
Calculate $\mu$ by using \eqref{linpromu};\\
Use additionally \eqref{alphaffpridua}, \eqref{muaff} and
\eqref{sigheu} to calculate $\sigma$; \\
Solve \eqref{linfindi} for $(\Delta x, \Delta \lambda, \Delta s)$;\\
Calculate $\alpha^{pri}$ and $\alpha^{dual}$ as in \eqref{alphpridua}; \\ 
Set $x^{k+1} = x^{k} + \alpha_{k}^{pri} \Delta x$;\\
Set $(\lambda^{k+1}, s^{k+1}) = (\lambda^{k}, s^{k}) +
\alpha_{k}^{dual}(\Delta \lambda, \Delta s)$;\\
$k = k+1$; \\
Set $r^{k} = (Ax^{k}-b, A^{\top}\lambda^{k} + s^{k} - c, XSe)$; \\
\textbf{until} $\text{norm}(r^{k}) < \epsilon$ (for a given $\epsilon
>0$).\\\hline 
\end{tabular}
\caption{Formal specification of a practical algorithm}
\label{lpppca}
\end{table}

As there are examples for that this algorithm diverges, no convergence theory
is available for the algorithm. Simple safeguards could be incorporated into
the method to force it into the convergence framework of existing methods or
to improve its robustness, but many practical codes do not implement these
safeguards because failures are rare. \label{IPlp}


\section{Extensions To Convex Quadratic Programming}\label{etcqp}

The interior point approach, introduced for linear programming in the previous
section, can also be applied to convex quadratic programs
through simple extensions of the linear programming algorithm. 

To keep the description of the interior point method simple, we consider a
QP with only inequality constraints:\footnote{If equality constraints are
  also present, they can be incorporated by simple extensions to the method
  described below}

\begin{subequations}\label{qpp}
\begin{align}
\min \frac{1}{2} x^{\top}Qx + d^{\top}x\\
\text{subject to}\;\; Ax \geq c ,
\end{align}
\end{subequations}

where $Q$ is  symmetric and positive definite, $d$ and $x$ are vectors in
$\mathbb R^{n}$, $c$ is a vector in $\mathbb  R^{m}$, and $A$ is an $m \times
n$ matrix.

\subsection{The KKT system and Newton's method}

The KKT conditions for \eqref{qpp} are :

\begin{subequations}\label{qppkkt}
\begin{align}
A^{\top}\lambda - Qx &= d, \\
Ax - s &= c, \\
\lambda_{i}s_{i} &= 0,\;\;\;\; i = 1,2,\ldots,m, \\
\lambda &\geq 0, \label{ipqpl}\\
s &\geq 0. \label{ipqps}
\end{align}
\end{subequations}

Since $Q$ is positive definite, these KKT conditions are necessary and
sufficient, and therefore we can solve \eqref{qpp} by finding solutions of
\eqref{qppkkt}.

Given a current iterate $(x,s,\lambda)$ that satisfies \eqref{ipqpl} and
\eqref{ipqps}, we can define the complementary measure $\mu$ as

\begin{align}
\mu = \frac{1}{m}\sum_{i=1}^{m}\lambda_{i}s_{i} =
\frac{\lambda^{\top}s}{m}. \label{qpmu} 
\end{align}

As in the previous section, we derive a practical, path-following, primal-dual
algorithm by considering the perturbed KKT conditions:

\begin{align}
F(x,\lambda,s; \sigma \mu) = \begin{bmatrix} A^{\top}\lambda - Qx - d \\ Ax - c
  - s \\ \Lambda S e - \sigma \mu e\end{bmatrix} &= 0 \label{ipqf}
\end{align}

where

\begin{align*}
\Lambda &= diag(\lambda_{1},\ldots,\lambda_{n}), \\
S &= diag(s_{1},\ldots,s_{n}),
\end{align*}

and $\sigma \in [0,1]$. The solutions of \eqref{ipqf} for all positive values
$\sigma$ and $\mu$ define the central path. This is a trajectory that leads to
the solution of the quadratic program as $\sigma \mu$ tends to zero.

By applying Newton's method to \eqref{ipqf}, we obtain the linear system

\begin{align}
\begin{bmatrix} Q & 0 & -A^{\top} \\ A & -I & 0 \\ 0 & \Lambda & S \end{bmatrix}
\begin{bmatrix} \Delta x \\ \Delta s \\ \Delta \lambda \end{bmatrix} = 
\begin{bmatrix} -r_{d} \\ -r_{c} \\ -\Lambda S e + \sigma \mu e \end{bmatrix} \label{nssm}
\end{align}

where

\begin{align*}
r_{c} &= Ax - s - c, \\ 
r_{d} &= Qx - A^{\top}\lambda + d.
\end{align*}

\subsection{Finding the optimal step length}

We define the new iterate as

\begin{subequations}\label{qpni}
\begin{align}
(x^{k+1}, s^{k+1}) &= (x^{k},s^{k}) + \alpha_{k}^{pri}(\Delta x^{k}, \Delta
s^{k}), \\ 
\lambda^{k+1} &= \lambda^{k} + \alpha_{k}^{dual}\Delta \lambda^{k}, 
\end{align}
\end{subequations}

where $(\alpha_{k}^{pri}, \alpha_{k}^{dual})$ are selected so as to
(approximately) minimize the optimality measure

\begin{subequations}\label{om}
\begin{align}
\|Qx^{k+1} - A^{\top}\lambda^{k+1} + d \|^{2}_{2} &+ \|Ax^{k+1} - s^{k+1} - c
\|^{2}_{2} + (s^{k+1})^{\top}\lambda^{k+1}, \\
\text{subject to} \;\; 0 &\leq \alpha_{k}^{pri} \leq \alpha_{k,\tau}^{pri}, \\
0 &\leq \alpha_{k}^{dual} \leq \alpha_{k,\tau}^{dual},  
\end{align}
\end{subequations}

where $x^{k+1},s^{k+1}$ and $\lambda^{k+1}$ are defined as functions of the step
lengths through \eqref{qpni} and $\alpha_{k,\tau}^{pri}$ and
$\alpha_{k,\tau}^{dual}$ are given by

\begin{align*}
\alpha_{k,\tau}^{pri} &= \max\{\alpha \in (0,1]:s^{k} + \alpha \Delta s^{k}
\geq (1 - \tau_{k})s^{k} \}, \\
\alpha_{k,\tau}^{dual} &= \max\{\alpha \in (0,1]:\lambda^{k} + \alpha \Delta
\lambda^{k} \geq (1 - \tau_{k})\lambda^{k} \},
\end{align*}

whereby the parameter $\tau_{k} \in (0,1)$ controls how far we back off from the
maximum step for which the conditions $s^{k} + \alpha \Delta s^{k} \geq 0$ and
$\lambda^{k} + \alpha \Delta \lambda^{k} \geq 0$ are satisfied.



\subsection{Choosing a starting point}

As for linear programming, the efficiency and
robustness of our practical algorithm can be greatly enhanced by choosing a
good starting point. Therefore we want to present at least a simple heuristic
that improves the choice of the starting point by moving an initial starting
point $(\overline{x}, \overline{s}, \overline{\lambda})$ form the user far
enough away from the boundary of the region $(s,\lambda) \geq 0$ to permit the
algorithm to take long steps on early iterations. Therefore, our heuristic
computes the affine scaling step $(\overline{x}^{aff}, \overline{s}^{aff},
\overline{\lambda}^{aff})$ from $(\overline{x},
\overline{s},\overline{\lambda})$ and  then sets

\begin{subequations}\label{qpsp}
\begin{align}
s_{0} &= \max(1, | \overline{s} + \Delta \overline{s}^{aff} |),\\
\lambda_{0} &= \max(1, | \overline{\lambda} + \Delta  \overline{\lambda}^{aff}
|), \\ 
x_{0} &= \overline{x}.
\end{align}
\end{subequations}

\subsection{A practical primal-dual algorithm}

The most popular practical algorithms for convex quadratic programming are, as
practical interior point methods for linear programming, based on Mehrotra's
predictor-corrector idea (for details see \cite{meh}). 

Therefore we first compute an affine scaling step $(\Delta
x^{aff}, \Delta \lambda^{aff}, \Delta s^{aff})$ by setting $\sigma = 0$ in
\eqref{nssm}. The following scheme calculates the maximum allowable
step lengths along the affine-scaling direction: 

\begin{subequations}\label{slaff}
\begin{align}
\alpha_{aff}^{pri} &\stackrel{\mathrm{def}}= \min(1, \min_{i:\Delta
  s_{i}^{aff}<0}- \frac {s_{i}}{\Delta s_{i}^{aff}}), \\
\alpha_{aff}^{dual} &\stackrel{\mathrm{def}}= \min(1, \min_{i:\Delta
  \lambda_{i}^{aff} < 0} - \frac {\lambda_{i}}{\Delta \lambda_{i}^{aff}}).
\end{align}
\end{subequations}

Using the above definitions, we set $\mu_{aff}$ in accordance with the
definition of $\mu$ in \eqref{qpmu} to be:

\begin{align}
\mu_{aff} = (s + \alpha_{aff}^{pri} \Delta s^{aff})^{\top}(\lambda +
\alpha_{aff}^{dual} \Delta \lambda^{aff})/n. \label{muaffqp}
\end{align}

The centering parameter $\sigma$ is chosen according to the following
heuristic:\footnote{\eqref{qpsigheu} has no solid analytical justification, but
  appears to work well in practice}

\begin{align}
\sigma = (\frac{\mu_{aff}}{\mu})^{3}. \label{qpsigheu}
\end{align}

Next we define the corrector step that aims to improve the affine scaling step
as 

\begin{align}
\begin{bmatrix} Q & 0 & -A^{\top} \\ A & -I & 0 \\ 0 & \Lambda & S \end{bmatrix}
\begin{bmatrix} \Delta x^{cor} \\ \Delta s^{cor} \\ \Delta \lambda^{cor} \end{bmatrix} = 
\begin{bmatrix} 0 \\ 0 \\ -\Delta X^{aff} \Delta S^{aff} e \end{bmatrix}. 
\end{align}

Finally, the total step is obtained by solving the following system:

\begin{align}
\begin{bmatrix} Q & 0 & -A^{\top} \\ A & -I & 0 \\ 0 & \Lambda & S \end{bmatrix}
\begin{bmatrix} \Delta x \\ \Delta s \\ \Delta \lambda \end{bmatrix} = 
\begin{bmatrix} -r_{c} \\ -r_{b} \\ -XSe - \Delta X^{aff} \Delta S^{aff} e +
  \sigma \mu e \end{bmatrix}. \label{qpts}
\end{align}

Finally we put together the different, above mentioned ingredients for a
practical successful algorithm in Table \ref{qpppca} below.\\

\begin{table}[h]
\small 
\begin{tabular}{l} \hline
\multicolumn{1}{c}{Practical Predictor-Corrector Algorithm} \\\hline
Calculate $(x^{0},\lambda^{0},s^{0})$ using an initial starting point from the
user $(\overline{x}, \overline{s}, \overline{\lambda})$ and \eqref{qpsp};
\\\hline 
k = 0;\\
\textbf{repeat}\\
Set $(x,\lambda,s)$ = $(x^{k},\lambda^{k},s^{k})$ and calculate $(\Delta
x^{aff}, \Delta \lambda^{aff}, \Delta s^{aff})$ by solving \eqref{nssm} with
$\sigma = 0$;\\ 
Calculate $\mu$ by using \eqref{qpmu};\\
Use additionally \eqref{slaff}, \eqref{muaffqp} and
\eqref{qpsigheu} to calculate $\sigma$;\\
Solve \eqref{qpts} for $(\Delta x, \Delta \lambda, \Delta s)$;\\
Select $\alpha_{k}^{pri}$ and $\alpha_{k}^{dual}$ to be the (approximate)
minimizers of the optimality measure in \eqref{om};\\ 
Set $x^{k+1} = x^{k} + \alpha_{k}^{pri} \Delta x$;\\
Set $(\lambda^{k+1}, s^{k+1}) = (\lambda^{k}, s^{k}) +
\alpha_{k}^{dual}(\Delta \lambda, \Delta s)$;\\
$k = k+1$;\\
Set $r^{k} = (Ax - s - c ,Qx - A^{\top}\lambda + d, \Lambda S e)$ \\
\textbf{until} $\text{norm}(r^{k}) < \epsilon$ (for a given $\epsilon
>0$).\\\hline
\end{tabular}
\caption{Formal specification of a practical algorithm}
\label{qpppca}
\end{table}


\chapter{Feasible Active-Set Methods}\label{cfasm}

This chapter deals with the description of practically successful active set
methods for convex quadratic programming. In the following sections we
consider the most important aspects of feasible active-set methods like the
working set, the subproblems at each iteration, the smart choice of a starting
point and the usage of updating factorizations. We conclude the chapter we a
comparison of active-set and interior point methods.

For the description of feasible active-set methods we use the following
problem formulation:\footnote{This formulation is equivalent to the other
  formulations of general QPs presented in this thesis. For details see
  Section \ref{dpf}.}

\begin{subequations}\label{pfas}
\begin{align}
\min_{x} \; \frac {1}{2}x^{\top}Qx &+ x^{\top}d \label{gqpof}\\
\text{subject to  } a_{i}^{\top}x &= c_{i}, \;\;\; i \in
\epsilon, \label{gqpec}\\ 
a_{i}^{\top}x &\geq c_{i}, \;\;\; i \in \iota, \label{gqpic}
\end{align}
\end{subequations}

where $Q$ is a symmetric, positive definite $n \times n$ matrix, $\epsilon$
and $\iota$ are finite sets of indices, and $d$, $x$ and $\{a_{i}\}, \;\; i \in
\epsilon \cup \iota ,$ are vectors in $\mathbb R^{n}$.

\section{Active Sets And Working Sets}

We now describe active-set methods for solving the quadratic program, given by
\eqref{pfas}.

If the contents of the optimal active set $A(x^{*})$, given by

\begin{align*}
A(x^{*}) = \{i \in \epsilon \cup \iota | a_{i}^{\top}x^{*} = c_{i}\}
\end{align*}

were known in advance, we could find the solution $x^{*}$ easily. Of course, we
usually do not have prior knowledge of $A(x^{*})$ and therefore determination of
this set is the main challenge facing active-set algorithms for quadratic
programs.

The simplex method starts by making a guess of the optimal active set, then
repeatedly uses gradient and Lagrange multiplier information to drop one index
from the current estimate of $A(x^{*})$ and add a new index, until optimality is
detected. Active-set methods for quadratic programs differ from the simplex
method in that the iterates and the solution $x^{*}$ are not necessarily
vertices of the feasible region.

There are primal, dual and primal-dual versions of active-set methods. We will
explain now primal methods, which are steadily decreasing the objective
function \eqref{gqpof} while remaining feasible with respect to the primal
problem.

Primal active-set methods find a step from one iterate to the next by solving
a quadratic subproblem in which some of the inequality constraints
\eqref{gqpic}, additionally to the equations \eqref{gqpec}, are treated as
equalities. This set of equations is called the working set and is denoted as
$W_{k}$ at the $k$th iterate $x_{k}$. We further assume that the gradients
$a_{i}$ of the constraints in $W_{k}$ are linearly independent.\footnote{If we
  use a linearly independent subset of the gradients $a_{i},i \in 
\{1,\ldots,m\}$ as initial working set $W_{0}$, the definition of the step
length in \eqref{defak} ensures that the linear independence is maintained for
the subsequent working sets $W_{k},k \geq 1$.}

\section{The Subproblems}\label{asmsp}

The first step in every iteration is to check whether the current iterate
$x_{k}$ minimizes \eqref{gqpof} in the subspace defined by $W_{k}$. If this is
not the case, we solve an equality-constrained quadratic subproblem, in which
the constraints belonging to $W_{k}$ are included and the other inequality
constraints are temporarily disregarded, to determine a step $p$, defined as

\begin{align}
p = x - x_{k}, \label{pdef}
\end{align}

Now, by substituting \eqref{pdef} in \eqref{gqpof}, we get

\begin{align*}
\min_{p} \; \frac {1}{2} p^{\top}Qp + g_{k}p + \phi_{k},
\end{align*}

where

\begin{align*}
g_{k} &= Qx_{k} + d, \\
\phi_{k} &= \frac {1}{2}x_{k}^{\top}Qx_{k} + d^{\top}x_{k},
\end{align*}

are independent of $p$. Therefore the subproblem to be solved at the $k$th
iteration can be written as

\begin{subequations}\label{pmin}
\begin{align}
\min_{p} \frac {1}{2} p^{\top}Qp + g_{k}p \label{spof}\\
\text{subject to} \;\; a_{i}^{\top}p = 0, \;\; i \in W_{k}. \label{spec}
\end{align}
\end{subequations}

We can solve the subproblem for example by a symmetric indefinite
factorization or by the Schur-complement method or by the Null-Space method
(for details see Section \ref{ecdm}).For the solution of this subproblem,
denoted by $p_{k}$, we have 

\begin{align}
a_{i}^{\top}(x_{k} + \alpha p_{k}) = a_{i}^{\top}x_{k} = c_{i}, \;\; \forall
\alpha, \;\; i \in W_{k}.
\end{align}

If the direction $p_{k}$ is nonzero, the objective function is strictly
decreasing (as $Q$ is positive definite) along this direction (for a proof see,
for example, Nocedal and Wright \cite[Theorem 16.6.]{nowr}). 

Now we have to decide how far to move along the direction $p_{k}$. We set 

\begin{align}
x_{k+1} = x_{k} + \alpha_{k} p_{k},
\end{align}

where we choose the step-length parameter $\alpha_{k}$ in order to maximize
the decrease in \eqref{gqpof} to be the largest value
in the range $[0,1]$ for which all constraints are satisfied:

\begin{align}
\alpha_{k} \stackrel{\mathrm{def}}= \min(1, \min_{i \notin W_{k},
  a_{i}^{\top}p_{k} < 0} \frac {c_{i} - a_{i}^{\top}x_{k}}{a_{i}^{\top}
  p_{k}}) \label{defak}
\end{align}

The constraint $i$, for which the minimum in \eqref{defak} is achieved, is
called blocking constraint. It is also possible for $\alpha_{k}$ to be zero,
because some constraint $i$ could fulfill $a_{i}^{\top}p_{k} < 0$  and
additionally be active at $x_{k}$ without belonging to $W_{k}$. 

If $\alpha_{k} < 1$, a new working set $W_{k+1}$ is constructed by adding one
of the blocking constraints to $W_{k}$. 

We continue to do this until we reach a point $\widehat{x}$ that minimizes
\eqref{pmin} over its current working set $\widehat{W}$. Such a point
$\widehat{x}$ satisfies the KKT conditions for the subproblem:

\begin{align*}
\begin{bmatrix} Q & A_{k}^{\top}\\ A_{k} &  0 \end{bmatrix}
\begin{bmatrix} -\widehat{p}\\ \widehat{\lambda}^{*} \end{bmatrix} =
\begin{bmatrix} g\\ h \end{bmatrix},
\end{align*}

where 

\begin{align*}
g &= d + Q\widehat{x}, \\
h &= A_{k}\widehat{x} - c,
\end{align*}

and $A_{k}$ is the Jacobian of the constraints in \eqref{pmin}
and $c_{k}$ is the vector whose components are $c_{i}, i \in W_{k}$. 
Furthermore $p = 0$ at $\widehat{x}$ and therefore we have that

\begin{align*}
\sum_{i \in \widehat{W}}a_{i}\widehat{\lambda_{i}} = g = Q \widehat{x} + d,
\end{align*}

for some Lagrange multipliers $\widehat{\lambda_{i}}, \;\; i \in
\widehat{W}$. It follows that $\widehat{x}$ and $\widehat{\lambda}$ satisfy
the first three KKT conditions for the original quadratic program \eqref{pfas}

\begin{subequations}
\begin{align}
Qx^{*} + c - \sum_{i \in A(x^{*})}\lambda_{i}^{*}a_{i} &= 0, \\
a_{i}^{\top}x^{*} &= c_{i}, \;\; \forall i \in A(x^{*}), \\
a_{i}^{\top}x^{*} &> c_{i}, \;\; \forall i \in \iota \setminus A(x^{*}),\\
\lambda_{i}^{*} &\geq 0, \forall i \in \iota \cap A(x^{*}), \label{gqpkktb}
\end{align}
\end{subequations}

if we define 

\begin{align*}
\widehat{\lambda_{i}} = 0 \;\; \forall i \in \iota \setminus \widehat{W},
\end{align*}

and consider the step length control defined in \eqref{defak}.

We now take a look at the fourth equation of the above KKT system
\eqref{gqpkktb}, which concerns the inequality constraints in
$\widehat{W}$. If these multipliers are all nonnegative, our solution
$(\widehat{x},\widehat{\lambda})$ is the global optimum for $\eqref{pfas}$. 

If, on the other hand, one or more multipliers are negative, the objective
function, given by \eqref{gqpof}, can be decreased by dropping one of these
constraints. Thus, the next step is to remove the most negative
multiplier\footnote{This choice is motivated by a sensitivity analysis
  concerning the removal of the Lagrange multipliers, which shows that the
  rate of decrease in the objective function is proportional to the negative
  magnitude of the multiplier. However the step length along the resulting
  direction may be small because of some blocking constraint. That's why the
  amount of decrease in the objective function is not guaranteed to be greater
than for other negative multipliers. Furthermore the magnitude of the
multipliers is dependent on the scaling of the corresponding
constraints. Therefore, as for the simplex method in linear programming,
strategies that are less sensitive to scaling often give better practical
results.} from $\widehat{W}$ and solve the subproblem, given by 
\eqref{pmin}, for the new working set.

It can be shown that the optimal value $p$ of this new subproblem gives a
direction that is feasible with respect to the dropped constraint (for a proof
see, for example, Nocedal and Wright \cite[Theorem 16.5.]{nowr}).

Hence, we have at least at every second iteration a direction $p_{k}$ that
guarantees together with the assumption that the step length $\alpha_{k}$ is
nonzero for every $p_{k} \not= 0$ that we have a strict decrease in the
objective function after two iterations. This fact finally guarantees finite
termination of our algorithm (for details see Nocedal and Wright \cite[Section
16.5]{nowr}).

\section{Choosing A Starting Point}

Various techniques can be used to determine an initial feasible point. One
such is to use a two-phase approach, where in Phase I an auxiliary linear
program is designed so that an initial basic feasible point is trivial to
find. This problem can be solved with the simplex method and its solution
gives a basic feasible point for the original \eqref{pfas} (for details see, for
example, Nocedal and Wright \cite[Section 13.5]{nowr}). 

An alternative approach is a penalty (or 'big $M$') method that includes a
measure of infeasibility in the objective function that is zero at the
solution. We introduce a scalar artificial variable $\eta$ into \eqref{pfas}
to get a measure of the constraint violation. So we solve the modified problem

\begin{subequations}\label{qpmp}
\begin{align}
\min_{x, \eta} \frac {1}{2}x^{\top}Qx + x^{\top}d + M\eta \\
\text{subject to  } (a_{i}^{\top}x - c_{i}) &\leq \eta, \;\;\; i \in
\epsilon\\
-(a_{i}^{\top}x - c_{i}) &\leq \eta, \;\;\; i \in
\epsilon\\
c_{i} - a_{i}^{\top}x &\leq \eta, \;\;\; i \in \iota \\
0 &\leq \eta,
\end{align}
\end{subequations}

for some large value of $M$. It can be shown by using the theory of Lagrange
multipliers (see Theorem \ref{conla} in Section \ref{conanal}) that if there
exist feasible points for the original problem \eqref{pfas}, then for $M$
sufficiently large, the solution of \eqref{qpmp} will have $\eta = 0$ and the
value of x will be also optimal for \eqref{pfas}.

To solve \eqref{pfas} we therefore use some heuristic to choose $M$, then solve
\eqref{qpmp} and increase $M$ if $\eta > 0$ until $\eta$
becomes zero. A feasible starting point for \eqref{qpmp}
can be obtained easily by just taking some guess $\tilde{x}$ and then choosing
$\eta$ large enough so that all constraints are satisfied. 

\section{Updating Factorizations}\label{udf}

In this subsection we explain an updating technique that is crucial to the
efficiency of the above presented active-set method.

As the working set can change by at most one index at
every iteration in the active-set method presented in this chapter, the KKT
matrix of the current iteration differs in at most one row and one column from
the KKT matrix of the previous iteration. Therefore we can compute the matrix
factors needed to solve the current subproblem by updating the factors
computed at the previous iteration. The total cost of the updating is in
general cheaper than solving the new system from the scratch.

We limit our discussion here to the null-space method, described in
\eqref{psp} - \eqref{lamopt}, but their are also ways to make an update for
the other methods presented in Subsection \ref{ecdm}. Suppose that the $m
\times n$ matrix $A$ has $m$ linearly independent rows and assume that the
orthogonal $n \times m$ 
matrix $Y$ and the orthogonal $n \times n-m$ matrix $Z$ are defined by means
of a QR factorization of $A^{\top}$ in the following way:\footnote{As $Z$ is
  not uniquely defined there are also other possible definitions of $Z$ (for
  details see Subsection \ref{ecdm})} 

\begin{align*}
A^{\top}\Pi = \begin{bmatrix}Y & Z \end{bmatrix} \begin{bmatrix} R \\
  0 \end{bmatrix}, 
\end{align*}

where $\Pi$ is a permutation matrix and $R$ is a square, nonsingular, upper
triangular $m \times m$ matrix.

Now let us take a look at the case where one constraint $a$ is added to the
working set. Our new constraint matrix $\overline{A}^{\top}$ has full column
rank and is equal to $[A^{\top} a]$. As $Y$ and $Z$ are orthogonal, we have

\begin{align}
\overline{A}^{\top} \begin{bmatrix} \Pi & 0 \\ 0 & 1 \end{bmatrix}
= \begin{bmatrix} A^{\top}\Pi & a \end{bmatrix} = \begin{bmatrix}Y &
        Z \end{bmatrix} \begin{bmatrix}R & S^{\top}a \\ 0 &
      \hat{Q}\begin{bmatrix} \gamma \\
        0 \end{bmatrix} \end{bmatrix}, \label{fms} 
\end{align}

where $\gamma$ is a scalar and $\hat{Q}$ is a orthogonal matrix that
transforms $Z^{\top}a$ in the following way:

\begin{align*}
\hat{Q}(Z^{\top}a) = \begin{bmatrix} \gamma \\ 0 \end{bmatrix}.
\end{align*}

From \eqref{fms} we can see that the new factorization has the form

\begin{align*}
\overline{A}^{\top} \overline{\Pi} = \begin{bmatrix} Y &
  Z\hat{Q}^{\top} \end{bmatrix} \begin{bmatrix} \overline{R} \\ 0 \end{bmatrix},
\end{align*}

where

\begin{align*}
\overline{\Pi} = \begin{bmatrix} \Pi & 0 \\ 0 & 1 \end{bmatrix}, \\
\overline{R} = \begin{bmatrix} R & Y^{\top}a \\ 0 & \gamma \end{bmatrix}.
\end{align*}

Now we choose $\overline{Z}$ to be the last $n - m - 1$ columns of
$Z\hat{Q}^{\top}$ to finish the update.

To update $Z$, we need to account for the cost of obtaining $\hat{Q}$ and the
cost for calculating $Z\hat{Q}$, which is of order $n(n-m)$. This is less
expensive than computing the new factors from scratch, which causes cost of
order $n^{2}m$, especially when the null space is small.

In the case that we want to remove an index from the working set, we have to
remove a row from $R$ and thus disturb its upper triangular property by
introducing a number of nonzeros on the diagonal immediately below the main
diagonal. We can restore the upper diagonal property by applying a sequence of
plane rotations that introduce a number of inexpensive transformations into
$Y$. The updated matrix $\overline{Z}$ is then the current matrix $Z$ augmented by
a single column $\overline{z}$:

\begin{align*}
\overline{Z} = \begin{bmatrix} \overline{z} & Z \end{bmatrix}.
\end{align*}

The total cost of the updating depend on the location of the removed column
but is in general cheaper than computing the QR factors from the scratch (for
details see Gill et al. \cite[Section 5]{gigo}).

Let us next consider the reduced Hessian $Z^{\top}QZ$. For problem
\eqref{pmin}, $h = 0$ in \eqref{ecqpkkt} and therefore $p_{y}$, given by
\eqref{defpy}, is also zero. Thus the equation for null-space vector $p_{z}$
reduces from \eqref{Zrs} to 

\begin{align*}
(Z^{\top}QZ)p_{z} = -Z^{\top}g. 
\end{align*} 

To update the Cholesky factorization of the reduced Hessian

\begin{align*}
Z^{\top}QZ = LL^{\top}
\end{align*}

a series of inexpensive, elementary operations can be used. Furthermore we can
update the reduced gradient $Z^{\top}g$ at the same time as $Z$(for details
see Nocedal and Wright \cite[Section 16.7]{nowr}).

\section{Comparison Of Active-Set And Interior Point Methods}

Interior point methods share common features that distinguish them from the
active set methods. Each interior point iteration is expensive to compute and
can make significant progress towards the solution, while the active set methods
usually require a large number of inexpensive iterations. Geometrically, the
active set methods for QP differ from the simplex method in that the iterates
are not necessarily vertices of the feasible region. Interior point methods
approach the boundary of the feasible set only in 
the limit. They may approach the solution either from the interior or exterior
of the feasible region, but they never actually lie on the boundary of this
region.

The numerical comparison of active-set and interior point methods for convex
quadratic programming, executed by Gould and Toint \cite{goto}, indicates that
interior-point methods are generally much faster on large problems. If some
warm start information is available, however, the active set methods are
generally preferable. Although a lot of research has been focused on improving
the warm-start ability of interior point methods, the full potential of
interior point methods in this area is not yet known.






\chapter{A Lagrangian Infeasible Active-Set Method}\label{iasm}

This chapter provides a description of a Lagrangian infeasible active-set
method for convex quadratic programming. In the following sections we
give detailed information about the different parts of the algorithm and then
take a look at the algorithm's convergence behaviour.

To describe the algorithm let $a$, $b$,
$d$ $\in$ $\mathbb R^n$, $c$ $\in$ $\mathbb 
R^m$, $A \in \mathbb R^{m \times n}$ and $Q$ $=$ $Q^T$ be given, with $Q$ a
positive definite $n$ $\times$ $n$ matrix. We consider a convex
quadratic minimization problem with equality constraints and simple bound
constraints:\footnote{This formulation is equivalent to the other
  formulations of general QPs presented in this thesis. For details see
  Section \ref{dpf}.}

\begin{subequations}\label{qpwbc}
\begin{align}
\min \; J(x) \;\; \text{subject to} \;\;
h(x) = 0 \;\; \text{and} \;\; b \leq x \leq a,
\end{align}
\end{subequations}

where

\begin{align*}
 J(x)\;&:=\;\frac{1}{2}x^TQx+d^Tx, \\
 h(x)\;&:=\;Bx-c.
\end{align*}

The KKT-system for \eqref{qpwbc} is given by

\begin{subequations}
\begin{align}
B^{\top}\lambda + Qx + d + s + t &= 0, \label{KKTfooc} \\
Bx &= c, \label{KKTec} \\
s \circ (x-b) &= 0, \\
t \circ (x-a) &= 0, \\
x-b &\geq 0, \\
a-x &\geq 0, \\
s &\leq 0, \\
t &\geq 0. 
\end{align}
\end{subequations}

It is well known that a vector $x$ together with vectors $\lambda$ $\in$
$\mathbb R^m$, $s$ $\in$ $\mathbb R^n$ and $t$ $\in$ $\mathbb R^n$ of Lagrange
multipliers for the equality and bound constraints furnishes a global
minimum of \eqref{qpwbc} if and only if ($x,\lambda,s,t$) satisfies the
KKT-system.

We now describe in some detail the approach sketched above. Therefore first we
give a survey of the main components of our algorithm in Table \ref{algsur}
and then take a closer look at the important parts of our approach in 
the following sections.\\

\begin{table}[h]
\small 
\begin{tabular}{l} \hline
\multicolumn{1}{c}{Prototype Algorithm} \\\hline
\textbf{Input:} Q symmetric, positive definite $n \times n$ matrix, A $n \times
m$ matrix, $a$, $b$, $d \in \mathbb R^{n}$,\\ \;\;\; $c \in \mathbb R^{m}$. $A_{1} \subseteq N$ and $A_{2}
\subseteq N$, $A_{1} \cap A_{2} = \emptyset$ e.g. $A_{1} = \emptyset$, $A_{2} = \emptyset$\\
\textbf{Output:} ($x, \lambda, s, t$) optimal solution \\\hline
\textbf{repeat until} ($x, \lambda, s, t$) is optimal \vspace{0.15cm}\\
\;\;\; Calculate the actual augmented Lagrange function.\vspace{0.15cm}\\

\;\;\; Minimize the augmented Lagrange function applying an infeasible active set
method.\\\;\;\; As initial active set use the optimal active set of the last
iteration.\vspace{0.15cm}\\

\;\;\; Update $\lambda$\vspace{0.15cm}\\

\;\;\; If there has been a change in the active set, try to solve\\ \;\;\;problem (P)
directly again using an infeasible active set method.\vspace{0.15cm}\\

\;\;\; Compare the solution of the direct approach with the one you have got\\ \;\;\;
from minimizing the augmented Lagrange function and take the ``better'' one. \\\hline
\end{tabular}
\caption{Description of the algorithm}
\label{algsur}
\end{table}

\section{Outer Algorithm: The Augmented Lagrangian Method}

We make use of the so-called augmented Lagrangian method, which was first
proposed by Hestenes \cite{hes} and Powell \cite{pow}, in our outer
algorithm. Therefore we define the augmented Lagrangian function, which is a
combination of the 
Lagrangian function and the quadratic penalty function, as:

\begin{align}
\mathscr{L}_{A}(x,\lambda;\sigma) = J(x) + \lambda^{\top}(Bx-c) + \frac{\sigma}{2}\|Bx-c\|^{2}  \label{defof}
\end{align}

We now try to solve the problem:

\begin{subequations}\label{iaprodef}
\begin{align}
 \min_{x} \; \mathscr{L}_{A}(x,\lambda;\sigma)\\
  \;\; \text{subject to} \;\; b \leq x \leq a ,
\end{align}
\end{subequations}

instead of the general quadratic program \eqref{qpwbc}.

Rewriting \eqref{iaprodef} by using \eqref{defof} gives

\begin{align*}
\min_{x} \; \frac{1}{2}x^{\top}\tilde{Q}x+\tilde{d}^{\top}x + e \\
\text{subject to} \;\; b \leq x \leq a
\end{align*}

where

\begin{align*}
\tilde{Q} = Q + \sigma B^{\top}B, \\
\tilde{d} = d + B^{\top}(\lambda+\sigma c),\\
e = \frac{\sigma}{2}c^{\top}c.
\end{align*}

Next we introduce an algorithm that fixes $\lambda$ at the current estimate
$\lambda_{k}$ at its $k$th iteration, fixes the penalty parameter
$\sigma$ to some well-chosen value, and performs minimization with respect to
$x$, of course considering the simple bound constraints for $x$. Using $x_{k}$
to denote the approximate minimizer of $\mathscr{L}_{A}(x,\lambda_{k};\sigma)
+ (s_{k} + t_{k})x$, we have by the first order optimality conditions that

\begin{align}
0 \approx Qx_{k} + d + B^{\top}(\lambda_{k}-\sigma (B^{\top}x_{k}-c)) + s_{k}
+ t_{k} 
\end{align}

Comparing this with the first order optimality condition for the general QP
\eqref{qpwbc}, given by \eqref{KKTfooc}, we get

\begin{align}
\lambda^{*} \approx \lambda_{k} - \sigma(Bx_{k}-c).
\end{align}

Therefore we update $\lambda$ by the rule

\begin{align}
\lambda_{k+1} = \lambda_{k} - \sigma(Bx_{k}-c) \label{lamud}
\end{align}

This first order updating rule is, for example, given by Nocedal and Wright
\cite{nowr} in formula 
(17.39) or by Bertsekas \cite{ber}. But they deduced it without considering
bound constraints.

We will show later on (in Theorem \eqref{conla}) that (under some conditions)
we can solve problem \eqref{qpwbc} by iteratively solving problem
\eqref{iaprodef} and updating $\lambda$.\\

\section{Inner Algorithm For Minimizing The Augmented Lagrange Function}

To solve problem \eqref{iaprodef} we use an infeasible active set method. This
method 
was already sucessfully applied to constrained optimal control problems (see
Bergounioux et al. \cite{beha, beit}) and to convex quadratic problems with
simple bound constraints (see \cite{kure}).

First we take a look at the KKT system for problem \eqref{iaprodef}:

\begin{subequations}\label{KKTia}
\begin{align}
 \tilde{Q}x + \tilde{d} + s + t &= 0 \label{LKKTM} \\
 s \circ (x-b) &= 0 \label{LKKTC1}\\
 t \circ (a-x) &= 0 \label{LKKTC2} \\
 x-b &\geq 0  \label{LKKTPC}\\
 a-x &\geq 0 \\
 s &\leq 0  \label{LKKTDC}\\
 t &\geq 0
\end{align}
\end{subequations}

 The crucial step in solving \eqref{iaprodef} is to identify those inequalities
 which are active on the lower bound and those which are active on the upper
 bound, i.e. the sets $A_{1} \subseteq N$ and $A_{2} \subseteq N$ ($A_{1} \cap
 A_{2} = \emptyset$), where the solution to \eqref{iaprodef} satisfies $x_{A_{1}}
 = b_{A_{1}}$ and $x_{A_{2}} = a_{A_{2}}$. Then, with $I := N \setminus (A_{1}
 \cup  A_{2})$, we must have $s_I = 0$, $t_{I} = 0$, $s_{A_{2}}=0$ and
 $t_{A_{1}}=0$.

 To compute the remaining elements $x_I$, $s_{A_{1}}$ and $t_{A_{2}}$ of $x$,
 $s$ and $t$, we use \eqref{LKKTM} and partition the equations and variables
 according to $A_{1}$, $A_{2}$ and $I$:

\begin{align}
\begin{pmatrix} \tilde{Q}_{A_{1}} & \tilde{Q}_{A_{1},A_{2}} & \tilde{Q}_{A_{1},I}  \\
  \tilde{Q}_{A_{2},A_{1}} & \tilde{Q}_{A_{2}} & \tilde{Q}_{A_{2},I} \\  \tilde{Q}_{I,A_{1}} &
  \tilde{Q}_{I,A_{2}} & \tilde{Q}_{I} \end{pmatrix}
\begin{pmatrix} x_{A_{1}} \\ x_{A_{2}} \\ x_I \end{pmatrix} +
\begin{pmatrix} \tilde{d}_{A_{1}} \\ \tilde{d}_{A_{2}} \\ \tilde{d}_I \end{pmatrix} +
\begin{pmatrix} s_{A_{1}} \\ s_{A_{2}} \\ s_I \end{pmatrix} +
\begin{pmatrix} t_{A_{1}} \\ t_{A_{2}} \\ t_{I} \end{pmatrix} = 0
\end{align}

The third set of equations can be solved for $x_I$, because $\tilde{Q}_I$ is by assumption positive definite:

\begin{align}
x_I = -\tilde{Q}_I^{-1}(\tilde{d}_I + \tilde{Q}_{I,A_{1}}b_{A_{1}} + \tilde{Q}_{I,A_{2}}a_{A_{2}}).
\end{align}

Substituting this into the first and second set of equations implies 

\begin{align}
s_{A_{1}} = -\tilde{d}_{A_{1}} - \tilde{Q}_{A_{1},N}x\\
t_{A_{2}} = -\tilde{d}_{A_{2}} - \tilde{Q}_{A_{2},N}x
\end{align}

 If our guesses for $A_{1}$ and $A_{2}$ would have been correct, then $x_I \geq
 b_I$, $s_{A_{1}} \leq 0$ and $t_{A_{2}} \geq 0$ would have to hold. Suppose this is
 not the case. Then we need to make a new 'guess' for $A_{1}$ and $A_{2}$, which we
 denote by $A_{1}^+$ and $A_{2}^{+}$. Let us first look at $s_{A_{1}}$. If $s_i < 0$, this
 confirms our previous guess $i \in A_{1}$, so we include $i$ also in
 $A_{1}^+$. Consider now $t_{A_{2}}$. If $t_i > 0$, this
 confirms our previous guess $i \in A_{2}$, so we include $i$ also in
 $A_{2}^+$. Let us finally look at $x_{I}$. If $x_i < b_i$ we set $x_i = b_i$ in
 the next iteration and hence we include $i$ in $A_{1}^+$. On the other hand if
 $x_{i} > a_{i}$ we set $x_{i} = a_{i}$ in the next iteration and therefore we
 include $i$ in $A_{2}^{+}$. Formally we arrive at

\begin{subequations}\label{nas}
\begin{align}
A_{1}^+ := \{i : x_i < b_i \; \text{or} \; s_i < 0\}\\
A_{2}^+ := \{i : x_i > a_i \; \text{or} \; t_i > 0\}.
\end{align}
\end{subequations}

So in each step of this iterative approach, we maintain the first order
optimality condition and the complementary constraints associated to
problem \eqref{iaprodef}, given by \eqref{LKKTM}, \eqref{LKKTC1} and
\eqref{LKKTC2}. As inital active sets we 
  take the empty sets in the first iteration and the 'best' (in terms of norm
  minimization of the equality constraints) active sets, we have found so far,
  for all consecutive iterations. The iterates of the algorithm are well
  defined, because in each step we get a unique solution for all $A_{1}
  \subseteq N$ and $A_{2} \subseteq N$, due to $\tilde{Q} \succ 0$.

\section{Inner Algorithm For Solving The Problem Directly}

After solving the quadratic program \eqref{iaprodef} for the actual
$\lambda_{k}$, we try to solve our general quadratic program with equality
constraints and simple bound constraints, given by \eqref{qpwbc},
directly by making use of the active sets that belong to the optimal value
of \eqref{iaprodef} as initial active sets.  

We again use the infeasible active set method described above. Solving the
system consisting of the first order optimality condition, given by
\eqref{KKTfooc}, and the equality constraint, given by \eqref{KKTec}, under
the additional constraints that $x_{A_{1}} = b_{A_{1}}$,
$x_{A_{2}}=a_{A_{2}}$, $s_I = 0$, $t_{I} = 0$, $s_{A_{2}}=0$, $t_{A_{1}}=0$
leads to 

\begin{align}
\begin{pmatrix} Q_I & B_{I,M}^{\top}  \\  Q_{A_{1},I} & B^{\top}_{A_{1},M}  \\
  Q_{A_{2},I} & B^{\top}_{A_{2},M} \\ B_{M,I} & 0 \end{pmatrix}
\begin{pmatrix} x_I \\ \lambda \end{pmatrix} =
\begin{pmatrix} -d_I - Q_{I,A_{1}}b_{A_{1}} - Q_{I,A_{2}}a_{A_{2}}\\  -d_{A_{1}} - s_{A_{1}} -
  Q_{A_{1}}b_{A_{1}} - Q_{A_{1},A_{2}}a_{A_{2}} \\ -d_{A_{2}} - t_{A_{2}} - Q_{A_{2},A_{1}}b_{A_{1}} -
  Q_{A_{2}}a_{A_{2}}  \\c - B_{M,A_{1}}b_{A_{1}} - B_{M,A_{2}}a_{A_{2}} \end{pmatrix}.\label{KKTKR}
\end{align}

Making $x_{I}$ explicit in the first set of equations of \eqref{KKTKR}

\begin{align}
x_{I}=Q_{I}^{-1}(-B^{\top}_{I,M}\lambda - d_{I} - Q_{I,A_{1}}b_{A_{1}} - Q_{I,A_{2}}a_{A_{2}})\label{dxla}
\end{align}

and using this in the fourth set of equations of \eqref{KKTKR} gives

\begin{multline}
B_{M,I}Q_{I}^{-1}B^{\top}_{I,M}\lambda = B_{M,A_{1}}b_{A_{1}} +
B_{M,A_{2}}a_{A_{2}} - \\c - B_{M,I}(Q_{I}^{-1}(d_{I} + Q_{I,A_{1}}b_{A_{1}} +
Q_{I,A_{2}}a_{A_{2}})).\\\label{dlamim} 
\end{multline}

If $B_{M,I}Q_{I}^{-1}B^{\top}_{I,M}$ is invertible, \eqref{dlamim} can be solved
for $\lambda$:

\begin{multline*}
\lambda = (B_{M,I}Q_{I}^{-1}B^{\top}_{I,M})^{-1}(B_{M,A_{1}}b_{A_{1}} +
B_{M,A_{2}}a_{A_{2}} - \\c - B_{M,I}(Q_{I}^{-1}(d_{I} + Q_{I,A_{1}}b_{A_{1}} +
Q_{I,A_{2}}a_{A_{2}}))). \\
\end{multline*}

By using $\lambda$ in \eqref{dxla}, we can calculate $x_{I}$. Finally making
use of $\lambda$ and $x_{I}$ in the second and third set of equations of \eqref{KKTKR}
yields $s_{A_{1}}$ and $t_{A_{2}}$:

\begin{align}
s_{A_{1}} = -d_{A_{1}} - Q_{A_{1}}b_{A_{1}} - Q_{A_{1},A_{2}}a_{A_{2}} - Q_{A_{1},I}x_{I} -  B^{\top}_{A_{1},M}\lambda \\
t_{A_{2}} = -d_{A_{2}} - Q_{A_{2}}a_{A_{2}} - Q_{A_{2},A_{1}}b_{A_{1}} - Q_{A_{2},I}x_{I} -  B^{\top}_{A_{2},M}\lambda
\end{align}

If our guess for $A_{1}$ and $A_{2}$ would have been correct, then $b_I \leq
x_{I} \leq a_{I}$, $s_{A_{1}} \leq 0$ and $t_{A_{2}} \geq 0$ would have to
hold. If this is not the case we arrive at a new active sets $A_{1}^{+}$ and
$A_{2}^{+}$, formally defined by 

\begin{align}
A_{1}^+ := \{i : x_i < b_i \; or \; s_i < 0\}\\
A_{2}^+ := \{i : x_i > a_i \; or \; t_i > 0\}.
\end{align}

If we cannot go on with the direct approach, because
$B_{M,I}Q_{I}^{-1}B^{\top}_{I,M}$ is not invertible for our current active sets
$A_{1}$ and $A_{2}$ or because we have reached a maximum number of
iterations $K < \infty$ 
we start a new outer iteration by calculating a new augmented Lagrange function.

Computational experience with our method indicates that
typically only few (most of the time only one) outer iterations
(multiplier-updates) and also only few (most of the time less than ten) inner
iterations (minimization of the Lagrange function and trying to solve
\eqref{qpwbc} directly) are required to reach the optimal solution.


To investigate the convergence behaviour of the algorithm we
look at convergence results for the augmented Lagrangian method and we examine
the convergence of our inner algorithms.

\section{Convergence Analysis Of The Augmented Lagrangian
  Method}\label{conanal} 

In this section we give a convergence result for the augmented Lagrangian
method and then take a closer look on an assumption made in Theorem
\ref{conla}.

The following result, given by Bertsekas \cite{ber} \cite{ber1}, gives
conditions under which there is a minimizer of
$\mathscr{L}_{A}(x,\lambda;\sigma)$ that lies close to $x^{*}$ and gives error
bounds both for $x_{k}$ and the updated multiplier estimate $\lambda^{k+1}$
obtained from solving the subproblem at iteration $k$.

\begin{theo}
Let $x^{*}$ be a strict local minimizer and a regular point of
(P). Furthermore let $\overline{\sigma}$ be a positive scalar such that
$\nabla_{xx}^{2}\mathscr{L}_{A}(x^{*},\lambda^{*};\overline{\sigma}) \succ 0$. Then there exist
positive scalars $\delta$, $\epsilon$, and $\kappa$ such that:
\begin{enumerate}
\item For all ($\lambda_{k}$, $\sigma$) in the set D $\subset R^{m+1}$ defined by 
\begin{align}
D = \{(\lambda_{k}, \sigma) \mid \| \lambda_{k} - \lambda^{*}\| < \delta\sigma,
\overline{\sigma} \leq \sigma \},         
\end{align}
the problem

\begin{align}
\min_{x} \; \mathscr{L}_{A}(x,\lambda_{k};\sigma) \;\; \text{subject to}
\;\; b \leq x \leq a, \| x - x^{*} \| \leq \epsilon
\end{align}

has a unique solution $x_{k}$. Moreover, we have

\begin{align}
\|x_{k} - x^{*} \| \leq \kappa \| \lambda_{k} - \lambda^{*} \|/\sigma.
\end{align}

\item For all ($\lambda_{k}$, $\sigma$) $\in$ D, we have 

\begin{align}
\|\lambda_{k+1} - \lambda^{*} \| \leq \kappa\|\lambda_{k}-\lambda^{*}\|/\sigma
\end{align} 

where $\lambda_{k+1}$ is given by \eqref{lamud}.

\item For all ($\lambda_{k}$, $\sigma$) $\in$ D, the matrix
  $\nabla_{xx}^{2}\mathscr{L}_{A}(x_{k},\lambda_{k};\sigma) \succ 0$.\\
\end{enumerate}
\label{conla}
\end{theo}

For a proof of the above Theorem see Bertsekas \cite{ber1, ber}.

Now we examine the positive definiteness assumption of the above theorem closer.

\begin{theo}
Let $x^{*}$ be a strict local minimizer and a regular point of
(P). Then

\begin{align}
\nabla_{xx}^{2}\mathscr{L}_{A}(x^{*},\lambda^{*};\sigma) \succ 0
\Leftrightarrow \sigma > \max\{-e_{1},\ldots,-e_{m}\}
\end{align}

where $e_{1},\ldots,e_{m}$ are the eigenvalues of $\{\nabla
h(x^{*})^{'}[\nabla^{2}_{xx}\mathscr{L}_{A}(x^{*},\lambda^{*};0)]^{-1}h(x^{*})\}^{-1}$. 
\end{theo}

For a proof of the above Theorem see again Bertsekas \cite{ber}.

Finally we adapt the above result to our problem structure.
For problem \eqref{iaprodef}

\begin{align*}
\{\nabla h(x^{*})^{'}[\nabla^{2}_{xx}\mathscr{L}_{A}(x^{*},\lambda^{*};0)]^{-1}h(x^{*})\}^{-1}
= (AQ^{-1}A^{'})^{-1}.
\end{align*} 

Therefore the assumption

\begin{align*}
\nabla_{xx}^{2}\mathscr{L}_{A}(x^{*},\lambda^{*};\overline{\sigma}) \succ 0
\end{align*}

in Theorem \ref{conla} is fulfilled for all $\overline{\sigma} > 0$.

\section{Convergence Analysis Of The Kunisch-Rendl Method}

In this section we generalize the proof idea used for the finite step
convergence result for the Kunisch-Rendl method with only upper bounds (for a
proof see Kunisch and Rendl \cite{kure}) to the case where we have lower and
upper bounds. The main aim of this section is to argue why the proof idea
does not work in this more general case any more.

\subsection{Index partition}

To investigate the behaviour of the algorithm, we look at two consecutive
iterations. Suppose that some iteration is carried out with the active sets
$A_{1}^{k} \subseteq N$ and $A_{2}^{k} \subseteq N$ ($A_{1}^{k} \cap A_{2}^{k} =
 \emptyset$) (for $k \geq 1$), yielding ($x^{k},s^{k},t^{k}$) as solution of
the KKT system \eqref{KKTia} for the current active sets. According to
\eqref{nas}, the new active sets are

\begin{align*}
A_{1}^{k+1} &:= \{i : x_i^{k} < b_i \;\; \text{or} \;\; s_i^{k} < 0\}\\
A_{2}^{k+1} &:= \{i : x_i^{k} > a_i \;\; \text{or} \;\; t_i^{k} > 0\}.
\end{align*}

Let $(x^{k+1},s^{k+1},t^{k+1})$ denote the solution of the KKT system
\eqref{KKTia} for the active sets $A_{1}^{k+1}$ and $A_{2}^{k+1}$. To avoid
too many superscripts, we write

\begin{align*}
(A, B, x, s, t) \;\; \text{for} \;\; (A_{1}^{k}, A_{2}^{k}, x^{k}, s^{k},
t^{k}) \;\; 
\text{and} \;\; (C, D, y, u, v) \;\; \text{for} \;\; (A_{1}^{k+1}, A_{2}^{k+1},
x^{k+1},s^{k+1},t^{k+1})
\end{align*}

Given $A$ and $B$, we have the set of inactive variables $I = N \setminus
(A \cup  B)$ and we find that $x, s, t, C, D, u, v$ are determined by

\begin{align*}
x_{A} &= b_{A}, \;\; x_{B} = a_{B}, \;\; s_{I} = s_{B} = 0, \;\; t_{I} = t_{A}
= 0, \;\; \tilde{Q}x + \tilde{d} + s + t = 0\\
C &= \{i : x_i < b_i \;\; \text{or} \;\; s_i < 0\} \;\; \text{and} \;\; D := \{i
: x_i > a_i \;\; \text{or} \;\; t_i > 0\}\\
J &= N \setminus (C \cup  D), \;\; y_{C} = b_{C}, \;\; y_{D} = a_{D}, \;\;
u_{J} = u_{D} = 0, \;\; v_{J} = v_{C} = 0, \;\; \tilde{Q}y + \tilde{d} + u + v
= 0 
\end{align*}

The following partition of $N$ into mutually disjoint subsets will be useful
in our analysis. We first partition $A$ into

\begin{align}
S := \{i \in A: s_{i} \geq 0 \}
\end{align}

and $A \setminus S$ and $B$ into

\begin{align}
T := \{i \in B: t_{i} \leq 0 \}
\end{align}

and $B \setminus T$. The set $I$ is partitioned into

\begin{align}
U &:= \{i \in I: x_{i} < b_{i} \} \\
V &:= \{i \in I: x_{i} > a_{i} \}
\end{align}

and $I \setminus (U \cup V)$. In Table \ref{ip} we summarize the relevant
information about $x, s, t, y, u, v$ for this partition. A nonspecified entry
indicates that the domain of the associated variable cannot be constrained.\\

\begin{table}[h]
\begin{center}
\begin{tabular}{| c || c | c | c | c | c | c |} \hline
  & $s$ & $t$ & $u$ & $v$ & $x$ & $y$ \\\hline \hline 
$S$ & $\geq$ 0 & = 0 & = 0 & = 0 & = b & \\\hline
$T$ & = 0 & $\leq$ 0 & = 0 & = 0 & = a & \\\hline
$A \setminus S$ & < 0 & = 0 &  & = 0 & = b & = b \\\hline
$B \setminus T$ & = 0 & > 0 & = 0 &  & = a & = a \\\hline
$U$ & = 0 & = 0 &  & = 0 & < b & = b\\\hline
$V$ & = 0 & = 0 & = 0 &  & > a & = a\\\hline
$I \setminus (U \cup V)$ & = 0 & = 0 & = 0 & = 0 & $\geq b \; \land \; \leq a$ &
\\\hline 
\end{tabular}
\end{center}
\caption{Partition of index set $N$}
\label{ip}
\end{table}

On the basis of the above table, we define the sets $K$ and $L$ that give the
indices of lower and upper infeasibility of $y$:

\begin{subequations}\label{dk}
\begin{align}
K_{1}&:=\{ i \in S: y_{i} < b_{i} \} \\
K_{2}&:=\{i \in T: y_{i} < b_{i} \} \\ 
K_{3}&:=\{i \in I \setminus (U \cup V): y_{i} < b_{i} \} \\
K &= K_{1} \cup K_{2} \cup K_{3} 
\end{align}
\end{subequations}

\begin{subequations}\label{dl}
\begin{align}
L_{1}&:=\{ i \in T: y_{i} > a_{i} \} \\
L_{2}&:=\{i \in S: y_{i} > a_{i} \} \\ 
L_{3}&:=\{i \in I \setminus (U \cup V): y_{i} > a_{i} \} \\
L &= L_{1} \cup L_{2} \cup L_{3} 
\end{align}
\end{subequations}

\subsection{The merit function}

Let us define our merit function as

\begin{align}
L_{c,d}(x,s,t) = \tilde{J}(x) + \frac {c}{2}\|g(x)\|^{2} + \frac {d}{2}
\|h(x)\|^{2} \label{lcd}
\end{align}

where

\begin{align}
\tilde{J}(x) &= x^{\top}\tilde{Q}x + \tilde{d}^{\top}x\\
g(x) &= \max(b-x,0)\\
h(x) &= \max(x-a,0)
\end{align}

In the remainder of this section we shall investigate the value of \eqref{lcd}
along the iterates of the algorithm:

\begin{align}
L_{c,d}(y,u,v) - L_{c,d}(x,s,t).
\end{align}

First, we consider the changes of the objective function during consecutive
iterations. One can not expect a monotone decrease of $\tilde{J}(x)$ as the
iterates may be infeasible.

\begin{lem} Let $(x,s,t)$, $(y,u,v)$, $U$ and $V$ be given as above and $W :=
  (U \;\cup \;V)$. Then, we have 

\begin{align}
\tilde{J}(y) - \tilde{J}(x) = \frac {1}{2}(y-x)^{\top}\begin{pmatrix}
  \tilde{Q}_{W} & 0 \\ 0 & -\tilde{Q}_{\overline{W}} \end{pmatrix}(y-x).
\end{align}

\end{lem}

\textbf{Proof.} We use the $Q$-inner product, $\langle a,b \rangle_{Q} :=
a^{\top}Q b$, with the associated norm $\|a\|^{2}_{Q} := \langle a,a
\rangle_{Q}$ and get

\begin{align}
\tilde{J}(y) - \tilde{J}(x) = \frac {1}{2}\|y\|_{\tilde{Q}}^{2} - \frac
{1}{2}\|x\|_{\tilde{Q}}^{2} + z^{\top}\tilde{d},\label{djj}
\end{align}

where $z = y - x$. Using the following identity

\begin{align}
\|a\|_{Q}^{2} - \|b\|_{Q}^{2} = 2 \langle a-b,a \rangle_{Q} - \|a -
b\|_{Q}^{2},
\end{align}

on the right hand side of \eqref{djj} we obtain

\begin{align}
\tilde{J}(y) - \tilde{J}(x) = - \frac {1}{2} z^{\top}\tilde{Q}z +
z^{\top}(\tilde{Q}y+\tilde{d}) 
\end{align}

Considering that $\tilde{Q}y + \tilde{d} = -u - v$, we get

\begin{align}
\tilde{J}(y) - \tilde{J}(x) = - \frac {1}{2} z^{\top}\tilde{Q}z -
z^{\top}(u+v).
\end{align}

Now $u_{i} = v_{i} = 0$ for $i \in S \cup T \cup (I \setminus W)$ and
$z_{i} = 0$ on $(A \setminus S) \cup (B \setminus T)$. Therefore $z^{\top}(u +
v) = \sum_{i \in W} z_{i}(u_{i}+v_{i})$. Furthermore $u + v - s - t =
-\tilde{Q}z$ and $u_{i} + v_{i} - s_{i} - t_{i} = u_{i} + v_{i}$ for $i \in
W$, and hence  

\begin{align}
-\sum_{i \in W}z_{i}(u_{i} + v_{i}) = \sum_{i \in W}
z_{i}(\tilde{Q}z)_{i} = z^{\top}\begin{pmatrix} \tilde{Q}_{W, N}\\
  0 \end{pmatrix} z.
\end{align}

Summarizing, we see that

\begin{align*}
\tilde{J}(y) - \tilde{J}(x) = - \frac {1}{2} z^{\top}\begin{pmatrix}
  \tilde{Q}_{W} & \tilde{Q}_{W,\overline{W}} \\
  \tilde{Q}_{\overline{W}, W} & \tilde{Q}_{\overline{W}} \end{pmatrix} z + z^{\top}\begin{pmatrix} \tilde{Q}_{W} &
  \frac {1}{2}\tilde{Q}_{W,\overline{W}} \\
  \frac {1}{2} \tilde{Q}_{\overline{W}, W} & 0 \end{pmatrix} z \;\;\;
 \boxed{}
\end{align*}

\begin{lem} Let (x,s,t), (y,u,v), U, V, K and L be given as above. Then we
  have

\begin{align}
\|g(y)\|^{2}  - \|g(x)\|^{2} = \sum_{i \in K}|y_{i} - b_{i}|^{2} - \sum_{i \in
U}|x_{i} - b_{i}|^{2}
\end{align}

as well as

\begin{align}
\|h(y)\|^{2}  - \|h(x)\|^{2} = \sum_{i \in L}|y_{i} - b_{i}|^{2} - \sum_{i \in
V}|x_{i} - b_{i}|^{2}
\end{align}

\end{lem}

\textbf{Proof.} The claim follows from the fact that $x$ is infeasible on the
lower bound precisely on $U$ and on the upper bound precisely
on $V$ (see Table \ref{ip}). Moreover, by the definition of the sets $K$ and
$L$ (see \eqref{dk} and \eqref{dl}), the variable $y$ is infeasible on
the lower bound on $K$ and on the upper bound on $L$. \;\;\; $\boxed{}$

In summary we have proved the following result.

\begin{pro} For every two consecutive triples (x,s,t) and (y,u,v) we have

\begin{multline}
L_{c,d}(y,u,v) - L_{c,d}(x,s,t) = \frac {1}{2}(y-x)^{\top}\begin{pmatrix}
  \tilde{Q}_{W} & 0 \\ 0 & -\tilde{Q}_{\overline{W}} \end{pmatrix}(y-x) +\\ \frac{c}{2} \sum_{i \in K}|y_{i} - b_{i}|^{2}
+ \frac{d}{2} \sum_{i \in L}|y_{i} - b_{i}|^{2} - \frac{c}{2} \sum_{i \in
  U}|x_{i} - b_{i}|^{2} - \frac{d}{2} \sum_{i \in V}|x_{i} -
b_{i}|^{2}. \;\;\; \boxed{}\\ 
\end{multline}
\label{lylx}
\end{pro}

Let us introduce $\mu := \lambda_{min}(\tilde{Q}) > 0$ as the smallest
eigenvalue of $\tilde{Q}$ and then formulate Proposition \ref{2lcd}.

\begin{pro} For every two consecutive triples (x,s,t) and (y,u,v) we have

\begin{multline}
2(L_{c,d}(y,u,v) - L_{c,d}(x,s,t)) = \|\tilde{Q}\|\|z_{W}\|^{2} -
\mu\|z_{\overline{W}}\|^{2} + \\c \|z_{K}\|^{2} + d\|z_{L}\|^{2} - c\|z_{U}\|^{2} -
d\|z_{V}\|^{2} \\
\end{multline}

\label{2lcd}
\end{pro}

\textbf{Proof.} We first note that for $i \in K$ we have $x_{i} \geq b_{i}$
and for $i \in L$ we have $x_{i} \leq a_{i}$. Hence $0 < b_{i} - y_{i} \leq
x_{i} - y_{i}$ for $i \in K$ and $0 < y_{i} - a_{i} \leq y_{i} - x_{i}$ for $i
\in L$, and therefore

\begin{align*}
\sum_{i \in K}(y_{i} - b_{i})^{2} &\leq \|z_{K} \|^{2} \\
\sum_{i \in K}(y_{i} - a_{i})^{2} &\leq \|z_{L} \|^{2} .
\end{align*}

Furthermore we have $y_{U} = b_{U}$ and $y_{V} = a_{V}$, and hence

\begin{align*}
\sum_{i \in U}(x_{i} - b_{i})^{2} &= \|z_{U} \|^{2} \\
\sum_{i \in V}(x_{i} - a_{i})^{2} &= \|z_{V} \|^{2} .
\end{align*}

Using Proposition \ref{lylx} we get

\begin{multline*}
2(L_{c,d}(y,u,v) - L_{c,d}(x,s,t)) = \|\tilde{Q}\|\|z_{W}\|^{2} -
\mu\|z_{\overline{W}}\|^{2} + \\c \|z_{K}\|^{2} + d\|z_{L}\|^{2} - c\|z_{U}\|^{2} -
d\|z_{V}\|^{2}   \;\;\;\boxed{} \\
\end{multline*}

\subsection{The need to bound $\|z_{K}\|$ and $\|z_{L}\|$}

The next goal should be to bound $\|z_{K}\|$ and $\|z_{L}\|$ in terms of
$\|z\|$. On $K_{1}$ and $L_{2}$ we have 

\begin{align*}
s_{K_{1}} \geq 0, s_{L_{2}} \geq 0, t_{K_{1}} = t_{L_{2}} = 0, u_{K_{1}} =
u_{L_{2}} = 0, v_{K_{1}} = v_{L_{2}} = 0,
\end{align*} 

and therefore 

\begin{align*}
(\tilde{Q}z)_{K_{1}} = s_{K_{1}} \geq 0 \;\; \text{and} \;\;
(\tilde{Q}z)_{L_{2}} = s_{L_{2}} \geq 0.
\end{align*} 

On $L_{1}$ and $K_{2}$ we have 

\begin{align*}
s_{L_{1}} = s_{K_{2}} = 0, t_{L_{1}} \leq 0, t_{K_{2}} \leq 0, u_{L_{1}} =
u_{K_{2}} = 0, v_{L_{1}} = v_{K_{2}} = 0,
\end{align*} 

and therefore 

\begin{align*}
(\tilde{Q}z)_{L_{1}} = t_{L_{1}} \leq 0 \;\; \text{and} \;\;(\tilde{Q}z)_{K_{2}}
= t_{K_{2}} \leq 0.
\end{align*} 

On $K_{3}$ and $L_{3}$ we have 

\begin{align*}
s_{K_{3}} = s_{L_{3}} = 0, t_{K_{3}} = t_{L_{3}} = 0, u_{K_{3}} = u_{L_{3}} =
0, v_{K_{3}} = v_{L_{3}} = 0, 
\end{align*} 

and thus 

\begin{align*}
(\tilde{Q}z)_{K_{3}} = (\tilde{Q}z)_{L_{3}} = 0.
\end{align*} 

It follows that 

\begin{align*}
(\tilde{Q}z)_{K} &= \tilde{Q}_{K}z_{K} + \tilde{Q}_{K,\overline{K}}z_{\overline{K}}
= \begin{pmatrix} s_{K1} \\ t_{K2} \\ 0 \end{pmatrix}\\
(\tilde{Q}z)_{L} &= \tilde{Q}_{L}z_{L} + \tilde{Q}_{L,\overline{L}}z_{\overline{L}}
= \begin{pmatrix} t_{L1} \\ s_{L2} \\ 0 \end{pmatrix}.
\end{align*} 

Taking the inner product  with $z_{K}$ and $z_{L}$ respectively, we obtain

\begin{align*}
z_{K}^{\top}(\tilde{Q}z)_{K} &= z_{K}^{\top}\tilde{Q}_{K}z_{K} +
z_{K}^{\top}\tilde{Q}_{K,\overline{K}}z_{\overline{K}} =
s_{K_{1}}^{\top}z_{K_{1}} + t_{K_{2}}^{\top}z_{K_{2}} \\
z_{L}^{\top}(\tilde{Q}z)_{L} &= z_{L}^{\top}\tilde{Q}_{L}z_{L} +
z_{L}^{\top}\tilde{Q}_{L,\overline{L}}z_{\overline{L}} =
t_{L_{1}}^{\top}z_{L_{1}} + s_{L_{2}}^{\top}z_{L_{2}}
\end{align*}

where

\begin{align*}
s_{K_{1}}^{\top}z_{K_{1}} &\leq 0 \;\; \text{but} \;\; t_{K_{2}}^{\top}z_{K_{2}}
\geq 0 \;\; \text{and}\\
t_{L_{1}}^{\top}z_{L_{1}} &\leq 0 \;\; \text{but} \;\;
s_{L_{2}}^{\top}z_{L_{2}} \geq 0.
\end{align*}

Thus we cannot derive the equations

\begin{subequations}\label{etd}
\begin{align}
z_{K}^{\top}\tilde{Q}_{K}z_{K} &\leq
z_{K}^{\top}\tilde{Q}_{K,\overline{K}}z_{\overline{K}} \;\; \text{or} \;\;
z_{K}^{\top}(\tilde{Q}z)_{K} \leq 0 \;\; \text{and}\\
z_{L}^{\top}\tilde{Q}_{L}z_{L} &\leq
z_{L}^{\top}\tilde{Q}_{L,\overline{L}}z_{\overline{L}} \;\; \text{or} \;\;
z_{L}^{\top}(\tilde{Q}z)_{L} \leq 0
\end{align}
\end{subequations}

that we would need to bound $\|z_{K}\|$ and $\|z_{L}\|$ in terms of
$\|z\|$. If we could derive the above equations \eqref{etd}, the rest of the
proof would be very similar to the one for only upper bounds in the paper of
Kunisch and Rendl \cite{kure}. We would have to set $c := d := \|\tilde{Q}\| +
\mu$ and define the conditions (C1) and (C2) slightly differently as

\begin{align*}
\text{condition (C1)} \hspace*{2cm} 2*\text{cond}(\tilde{Q}) <  
(\frac {\mu}{\nu})^{2} - 2 \\ 
\text{condition (C2)} \hspace*{2cm} 2*\text{cond}(\tilde{Q}) <
(\frac {q}{r})^{2} - 2,
\end{align*}

where the diagonal matrix $D := diag(q_{11},\ldots,q_{nn})$ is consisting of
the main diagonal elements of $\tilde{Q}$ and 

\begin{align*}
r &:= \|Q-D\|,\\
\text{cond}(\tilde{Q}) &= \frac
{\lambda_{\max}(\tilde{Q})}{\lambda_{\min}(\tilde{Q})}, \\ 
\nu &:= \max\{\|\tilde{Q}_{A,\overline{A}}\| : A \subset N, A \not = 0, A
\not = N \}, \\
q &:= \min\{q_{ii} : i \in N\}.
\end{align*}

Although we cannot prove the convergence of the Lagrangian infeasible
active-set method presented in Chapter
\ref{iasm}, the method converges very fast in practice, as we will see in the
following section.

\section{Computational Experience}\label{scoex}

In this  we look at the practical behaviour of our algorithm
by considering a variety of test problems.

The only nontrivial inputs to our algorithm are the initial active sets
$A_{1}$ and $A_{2}$,
the initial Lagrange multiplier $\lambda_{0}$ and the penalty parameter
$\sigma$. Our algorithm is quite insensitive to their choice.\footnote{$A_{2}$
  and $A_{2}$ are
  chosen as empty sets, $\lambda_{0}$ as zero vector and $\sigma$ equal
  to 10000}

\subsection{Randomly Generated Dense Problems}

At first we study in some detail randomly generated problems, where $Q$ and
$B$ are dense matrices. We vary the number of variables $n$ and the number of
equality constraints $m$. In order for the reader to be able to reproduce some
of the following results, we provide the MATLAB commands that we used to
generate the data $Q$, $B$, $d$, $a$, $b$ and $c$.

>> n = 500; (or n = 1000 or ... or n = 15000) \\
>> m = 5; (in gernal: m = n/100; m = n/10; m = n/2) \\
>> rand('seed',n+m) \\
>> x = rand(n,1); \\
>> B = rand(m,n); \\
>> c = B*x; \\
>> d = rand(n,1); \\
>> Z = rand(n)-0.5; \\
>> Q = Z'*Z + eye(n); \\
>> b = zeros(n,1); \\
>> a = ones(n,1);\\

In Table \ref{ddt} below we summarize the key performance
features of our algorithm for different problem sizes:

\begin{itemize}
\item The number of $\lambda$-Updates, which is equal with the number of outer
  iterations,
\item the number of iterations we run to solve the inner problem
  \eqref{iaprodef},
\item the number of iterations we carry out to try to solve the underlying
  problem \eqref{qpwbc} directly,
\item the time that is needed on a workstation with 3 GHz and 10 GB RAM get the
  optimal solution.\\
\end{itemize}

\hspace*{1cm}

\begin{table}[h]
\small 
\begin{tabular}{|c|c|c|c|c|c|} \hline \hline
&  & $\lambda$-updates & Iter. inner problem & Iter. underlying problem &
seconds\\\hline
\multirow{2}{*}{ n=500 } & m=50 & 1(0) & 7(0.67) & 1(0) &  0.09(0.01)   \\\cline{2-6}
& m=250 & 1(0) & 7.9(0.74) & 1(0) & 0.2(0.02)\\\hline
\multirow{2}{*}{ n=1000} & m=100 & 1(0) & 7.9(0.57) & 1.1(0.32) & 0.5(0.06)\\\cline{2-6}
& m=500 & 1(0) & 8.7(0.82) & 1(0) & 1.0(0.04) \\\hline
\multirow{2}{*}{ n=3000} & m=300 & 1(0) & 9(0.47) & 1.1(0.32) & 9.7(1.15) \\\cline{2-6}
& m=1500 & 1(0)  & 9.5(0.71) & 1.6(0.52) & 23.9(4,4) \\\hline
\multirow{2}{*}{n=5000} & m=500 & 1 & 9 & 1 & 37.2\\\cline{2-6}
& m=2500 & 1 & 11 & 2 & 118 \\\hline
\multirow{2}{*}{n=10000} & m=1000 & 1 & 10 & 1 & 282 \\\cline{2-6}
& m=5000 & 1 & 10 & 2 & 864 \\\hline
\multirow{2}{*}{n=15000} & m=1500 & 1 & 10 & 1 & 963 \\\cline{2-6}
& m=7500 & 1 & 10 & 2 & 2896 \\\hline
\end{tabular}
\caption{Dense data: Key performance features of our algorithm for different
  problem sizes and structures}
\label{ddt}
\end{table}

\hspace*{1cm}

Until n=3000 we perform 10 runs for different random data and give the
expectation value and in parenthesis the standard deviation. For larger
problems we only make one run in order to save time.\footnote{The low
standard deviations justify this action.}

We can see that the algorithm always needs only one outer iteration. It takes
no more than 11 iterations to solve the inner problem \eqref{iaprodef} and no
more than 2 
further iterations to finally get the exact numerical solution for the
underlying problem \eqref{qpwbc}. If we compare these results with the ones we
obtained using the same algorithm on the same problem data but without
considering upper bounds, we recognize that we need more iterations to solve
the inner problem \eqref{iaprodef}, but less time, as the systems of equations
we have to solve are smaller (for further details see \cite[Section 4.1]{hun}).

\subsection{Randomly Generated Sparse Problems}

Next we study randomly generated problems, where $Q$ and
$B$ are sparse matrices with 10 nonzero entries per row in average. We vary
the number of variables $n$ and the number of equality constraints $m$. In order for the reader to be able to reproduce some
of the following results, we provide the MATLAB commands that we used to
generate the data $Q$, $A$, $d$, $b$ and $c$, where by $nz$ we denote the
average number of nonzero entries per row in $Q$ and $A$.

>> n = 500; (or n = 1000 or ... or n = 20000) \\
>> m = 5; (in gernal: m = n/100; m = n/10; m = n/2) \\
>> rand('seed',n+m) \\
>> x = rand(n,1);\\
>> B = rand(m,n);\\
>> for i = 1:m;\\
\hspace*{0.6cm} y = rand(n,1);\\
\hspace*{0.6cm} for j = 1:n;\\
\hspace*{1cm} if y(j) > nz/n;\\
\hspace*{1.4cm}  B(i,j) = 0;\\
\hspace*{1cm} end;\\
\hspace*{0.6cm} end;\\
\hspace*{0.2cm} end;\\
>> c = B*x;\\
>> d = rand(n,1);\\
>> Z = sprand(n,n,0.1);\\
>> Q = Z'*Z + eye(n);\\
>> for i = 1:n;\\
\hspace*{0.6cm} y = rand(n,1);\\
\hspace*{0.6cm} for j = i:n;\\
\hspace*{1cm} if y(j) > nz/n;\\
\hspace*{1.4cm} Q(i,j) = 0;\\
\hspace*{1.4cm} Q(j,i) = 0;\\
\hspace*{1cm} end;\\
\hspace*{0.6cm} end;\\
\hspace*{0.2cm} end;\\
>> Q = Q + (abs(min(eig(Q)))+1)*eye(n);\\
>> b = zeros(n,1);\\
>> a = ones(n,1);\\

In Table \ref{sdt} below we summarize the key performance features of our
algorithm for different problem sizes.\\

\hspace*{1cm}

\begin{table}[h]
\small 
\begin{tabular}{|c|c|c|c|c|c|} \hline \hline
&  & $\lambda$-updates & Iter. inner problem & Iter. underlying problem &
seconds\\\hline
\multirow{2}{*}{ n=500 } & m=50 & 1(0) & 4.3(0.48) & 1(0) & 0.06(0) \\\cline{2-6}
& m=250 & 1(0) & 5.6(0.70) & 1.7(0.48) & 0.2(0.04)\\\hline
\multirow{2}{*}{ n=1000} & m=100 & 1(0) & 4.8(0.42) & 1.3(0.48) & 0.3(0.04) \\\cline{2-6}
& m=500 & 1(0) & 6.2(0.79) & 2.1(0.57) & 1.3(0.22) \\\hline
\multirow{2}{*}{ n=3000} & m=300 & 1(0) & 5.7(0.67) & 1.9(0.32) & 6.3(0.55)\\\cline{2-6}
& m=1500 & 3.1(6.64) & 9.6(8.22) & 2.9(0.99) & 32.6(12.8) \\\hline
\multirow{2}{*}{n=5000} & m=500 & 1 & 6 & 2 & 25.1 \\\cline{2-6}
& m=2500 & 1 & 8 & 3 & 135 \\\hline
\multirow{2}{*}{n=10000} & m=1000 & 1 & 6 & 2 & 186 \\\cline{2-6}
& m=5000 & 1 & 6 & 4 & 1235 \\\hline
\multirow{2}{*}{n=15000} & m=1500 & 1 & 6 & 3 & 797 \\\cline{2-6}
& m=7500 & 1 & 8 & 4 & 4208 \\\hline
\end{tabular}
\caption{Sparse data: Key performance features of our algorithm for different
  problem sizes and structures}
\label{sdt}
\end{table}

\hspace*{1cm}

Until n=3000 we perform 10 runs for different random data and give the
expectation value and in parenthesis the standard deviation. For larger
problems we only make one run in order to save time.\footnote{The low
standard deviations justify this action.}

We can see that the algorithm most of the time needs only one outer
iteration.It takes
no more than 8 iterations to solve the inner problem \eqref{iaprodef} and no
more than 4 further iterations to finally get the exact numerical solution for
the underlying problem \eqref{qpwbc}. If we compare these results with the
ones we obtained using the same algorithm on the same problem data but without
considering upper bounds, we recognize that we need more iterations to solve
the inner problem \eqref{iaprodef}, but less time, as the systems of equations
we have to solve are smaller (for further details see \cite[Section 4.2]{hun}).


\chapter{Discussion}\label{cdis}

The main interest of this diploma thesis was to describe and compare different,
practically successful solution methods for general convex quadratic problems
with arbitrary linear constraints. 

Therefore we showed the equivalence
of different QP problem formulations and presented some important so-called
direct methods for solving equality-constrained QPs in Chapter
\ref{prelim}. After this, we covered the most
important aspects for practically successful interior point and
active-set methods for convex quadratic programming in Chapter \ref{cipm} and
Chapter \ref{cfasm} respectively. 

Finally, as the core of the diploma thesis, we presented a combination of the
augmented Lagrangian method with an infeasible active set method as a new
algorithm for solving QPs efficiently in Chapter \ref{iasm}.

Among the special features of this algorithm are its ability to find
the exact numerical solution of the problem and the fact that at each
iteration level the size of the linear system that must be solved is
determined by the currently inactive set that can be significantly smaller
than the total set of variables. As a consequence the proposed algorithm
differs significantly from the interior point methods that we described in
Chapter \ref{cipm}. Because of its ability to 
'correct' many active variables to inactive ones and vice versa in each
iteration and its computationally cheap definition of the new active sets, the
algorithm also seems preferable to the feasible active set 
methods presented in Chapter \ref{cfasm}.

From the numerical experiments in Section \ref{scoex}
we observe that the algorithm can mostly find the optimal solution in the
first try to solve the system directly. This is 
certainly one of its distinguishing practical features. Furthermore the total
number of 
iterations is frequently quite insensitive with respect to data and
initialization. 

The next step of research will be to compare our algorithm
with other fast software for solving QPs on different test problems.
 

\bibliographystyle{acm}
\bibliography{Dip_lit}

\begin{thebibliography}{10}

\bibitem{ilo}
{\sc 8.0, I.~C.}
\newblock {\em User's Manual}.
\newblock France, 2002.

\bibitem{anan}
{\sc Andersen, E., and Andersen, K.}
\newblock The mosek interior point optimizer for linear programming: an
  implementation of the homogeneous algorithm.
\newblock In {\em High Perfmormance Optimization}, T.~Frenk, K.~Roos, and
  S.~Zhang, Eds. Kluwer Academic Publishers, 2000, pp.~197--232.

\bibitem{beha}
{\sc Bergounioux, M., Haddou, M., Hinterm\"uller, M., and Kunisch, K.}
\newblock A comparison of interior point methods and a {M}oreau-{Y}osida based
  active set strategy for constrained optimal control problems.
\newblock {\em SIAM Journal on Optimization 11}, 2 (2000), 495--521.

\bibitem{beit}
{\sc Bergounioux, M., Ito, K., and Kunisch, K.}
\newblock Primal-{D}ual {S}trategy for {C}onstrained {O}ptimal {C}ontrol
  {P}roblems.
\newblock {\em SIAM Journal on Control and Optimization 37\/} (1999),
  1176--1194.

\bibitem{ber}
{\sc Bertsekas, D.}
\newblock {\em Constrained Optimization and Lagrange Multiplier Methods}.
\newblock Athena Scientific, Massachsetts, 1996.

\bibitem{ber1}
{\sc Bertsekas, D.}
\newblock {\em Nonlinear Programming}, 2~ed.
\newblock Athena Scientific, Massachsetts, 1999.

\bibitem{bumo}
{\sc Burke, V., and Mor\'{e}, J.}
\newblock Exposing constraints.
\newblock {\em SIAM Journal on Optimization 4\/} (1994), 573--595.

\bibitem{byno}
{\sc Byrd, R., Nocedal, J., and Waltz, R.}
\newblock Knitro: An integrated package for nonlinear optimization.
\newblock In {\em Large-Scale Nonlinear Optimization}, G.~Di~Pillo and M.~Roma,
  Eds. Springer, 2006.

\bibitem{cogo1}
{\sc Conn, A., Gould, N., and Toint, P.}
\newblock Testing a class of algorithms for solving minimization problems with
  simple bounds on the variables.
\newblock {\em Mathematics of Computation 50\/} (1988), 399--430.

\bibitem{cogo}
{\sc Conn, A., Gould, N., and Toint, P.}
\newblock {\em Lancelot: A Fortran Package for Large-scale Nonlinear
  Optimization (release A)}.
\newblock Springer, Heidelberg, New York, 1992.

\bibitem{domo1}
{\sc Dolan, E., and Mor\'{e}, J.}
\newblock Benchmarking optimization software with performance profiles.
\newblock {\em Mathematical Programming 91\/} (2002), 201--213.

\bibitem{domo}
{\sc Dolan, E., Mor\'{e}, J., and Munson, T.}
\newblock {\em Benchmarking Optimization Software with COPS 3.0}.
\newblock Technical Report ANL/MCS-TM-273, 2004.

\bibitem{dos}
{\sc Dost\'{a}l, Z.}
\newblock An optimal algorithm for bound and equality constrained quadratic
  programming problems with bounded spectrum.
\newblock {\em Computing 78\/} (2006), 311--328.

\bibitem{fle}
{\sc Fletcher, R.}
\newblock Stable reduced hessian updates for indefinite quadratic programming.
\newblock {\em Mathematical Programming 87\/} (2000), 251--264.

\bibitem{fogi}
{\sc Forsgren, A., and Gill, P.}
\newblock Primal-dual interior methods for nonconvex nonlinear programming.
\newblock {\em SIAM Journal on Optimization 8\/} (1998), 1132--1152.

\bibitem{fogi1}
{\sc Forsgren, A., Gill, P., and Wright, M.}
\newblock Interior methods for nonlinear optimization.
\newblock {\em SIAM review 44\/} (2003), 525 -- 597.

\bibitem{frmi}
{\sc Freund, R., and Mizuno, S.}
\newblock Interior point methods: Current status and future directions.
\newblock {\em Optima 51\/} (1996), 1--9.

\bibitem{gewr1}
{\sc Gertz, M., and Wright, S.}
\newblock {\em OOQP User Guide}.
\newblock Technical Memorandum No. 2520, 2001.

\bibitem{gewr}
{\sc Gertz, M., and Wright, S.}
\newblock Object-oriented software for quadratic programming.
\newblock {\em ACM Transactions on Mathematical Software (TOMS) 29\/} (2003),
  58--81.

\bibitem{gigo}
{\sc Gill, P., Golub, G., Murray, W., and Saunders, M.}
\newblock Methods for modifying matrix factorizations.
\newblock {\em Mathematics of Computation 28\/} (1974), 505--535.

\bibitem{gimu1}
{\sc Gill, P., and Murray, W.}
\newblock Numerically stable methods for quadratic programming.
\newblock {\em Mathematical Programming 14\/} (1978), 349--372.

\bibitem{gou}
{\sc Gould, N.}
\newblock On practical conditions for the existence and uniqueness of solutions
  to the general equality quadratic programming problem.
\newblock {\em Mathematical Programming 32\/} (1985), 90 -- 99.

\bibitem{goto1}
{\sc Gould, N., and Toint, P.}
\newblock An iterative working-set method for large-scale non-convex quadratic
  programming.
\newblock {\em Applied Numerical Mathematics 43\/} (2002), 109--128.

\bibitem{goto}
{\sc Gould, N., and Toint, P.}
\newblock Numerical methods for large-scale non-convex quadratic programming.
\newblock In {\em Trends in Industrial and Applied Mathematics}, A.~Siddiqi and
  M.~Ko\v{c}vara, Eds. Kluwer Academic Publishers, Dordrecht, 2002,
  pp.~149--179.

\bibitem{hes}
{\sc Hestenes, M.}
\newblock Multiplier and gradient methods.
\newblock {\em Journal of Optimization Theory and Applications 4\/} (1969),
  303--320.

\bibitem{hun}
{\sc Hungerl\"{a}nder, P.}
\newblock {\em A solution method for convex quadratic problems}.
\newblock Austria, 2008.

\bibitem{kar}
{\sc Karmarkar, N.}
\newblock A new polynomial-time algorithm for linear programming.
\newblock {\em Combinatorica 4\/} (1984), 373--395.

\bibitem{kha}
{\sc Khachiyan, L.}
\newblock A polynomial algorithm in linear programming.
\newblock {\em Soviet Mathematics Doklady 20\/} (1979), 191--194.

\bibitem{klmi}
{\sc Klee, V., and Minty, G.}
\newblock How good is the simplex algorithm?
\newblock In {\em Inequalities}, O.~Shisha, Ed. Academic Press, New York, 1972,
  pp.~159--175.

\bibitem{kure}
{\sc Kunisch, K., and Rendl, F.}
\newblock An infeasible active set method for convex problems with simple
  bounds.
\newblock {\em SIAM Journal on Optimization 14}, 1 (2003), 35--52.

\bibitem{meg}
{\sc Megiddo, N.}
\newblock Pathways to the optimal set in linear programming.
\newblock In {\em Progress in Mathematical Programming: Interior-Point and
  Related Methods}, N.~Megiddo, Ed. Springer, New York, 1989, pp.~131--158.

\bibitem{meh}
{\sc Mehrotra}.
\newblock On the implementation of a primal-dual interior point method.
\newblock {\em SIAM Journal on Optimization 2\/} (1992), 575--601.

\bibitem{mes}
{\sc M\'{e}sz\'{a}ros, C.}
\newblock The bpmpd interior point solver for convex quadratic problems.
\newblock {\em Optimization Methods and Software 11\/} (1999), 431--449.

\bibitem{mit}
{\sc Mittelmann, H.}
\newblock Benchmarking interior point lp/qp solvers.
\newblock {\em Optimization Methods and Software 11\/} (1999), 655--670.

\bibitem{musa}
{\sc Murtagh, B., and Saunders, M.}
\newblock Large-scale linearly constrained optimization.
\newblock {\em Mathematical Programming 14\/} (1978), 41--72.

\bibitem{musa1}
{\sc Murtagh, B., and Saunders, M.}
\newblock {\em MINOS 5.5 User's guide}.
\newblock Technical Report SOL 83-20R, Standford University, 1998.

\bibitem{nene}
{\sc Nesterov, Y., and Nemirovskii, A.}
\newblock {\em Interior Point Polynomial Methods in Convex Programming: Theory
  and Applications}.
\newblock SIAM, Philadelphia, PA, 1994.

\bibitem{neto}
{\sc Nesterov, Y., and Todd, M.}
\newblock Self-scaled barriers and interior point methods for convex
  programming.
\newblock {\em Mathematics of Operations Research 22\/} (1997), 1--42.

\bibitem{neto1}
{\sc Nesterov, Y., and Todd, M.}
\newblock Primal-dual interior-point methods for self-scaled cones.
\newblock {\em SIAM Journal on Optimization 8\/} (1998), 324--362.

\bibitem{nowr}
{\sc Nocedal, J., and Wright, S.}
\newblock {\em Numerical {O}ptimization}, 2~ed.
\newblock Springer, New York, 2006.

\bibitem{powr}
{\sc Potra, F., and Wright, S.}
\newblock Interior-point methods.
\newblock {\em Journal of Computational and Applied Mathematics 124\/} (2000),
  281--302.

\bibitem{pow}
{\sc Powell, M.}
\newblock A method for nonlinear constraints in minimization problems.
\newblock In {\em Optimization}, R.~Fletcher, Ed. Academic Press, New York,
  1969, pp.~283--298.

\bibitem{van1}
{\sc Vanderbei, R.}
\newblock Loqo: An interior point code for quadratic programming.
\newblock {\em Optimization Methods and Software 11\/} (1999), 451--484.

\bibitem{van}
{\sc Vanderbei, R.}
\newblock {\em Linear Programming: Foundations and Extension}, 2~ed.
\newblock Kluwer Academic, Norwell, MA, 2001.

\bibitem{vaum}
{\sc Vanderbei, R.}
\newblock {\em LOQO User's Manual - Version 4.05}.
\newblock Technical Report No. ORFE-99, Princeton University, 2006.

\bibitem{wapl}
{\sc Waltz, R., and Plantenga, T.}
\newblock {\em Knitro User's Manual: Version 5.1}.
\newblock 2007.

\bibitem{wri}
{\sc Wright, S.}
\newblock {\em Primal-Dual Interior-Point Methods}.
\newblock siam, Philadelphia, PA, 1997.

\bibitem{ye}
{\sc Ye, Y.}
\newblock {\em Interior Point Algorithms: Theory and Analysis}.
\newblock John Wiley and Sons, New York, 1997.

\end{thebibliography}

\end{document}